%
\input amstex 
\documentstyle{amsppt}
\input bull-ppt
\keyedby{bull389e/mhm}

\define\np{\noindent}
\define\ID{\operatorname{id}}
\define\GL{\operatorname{GL}}
\define\TR{\operatorname{tr}}
\define \BK{\bold{K}}


\topmatter
\cvol{28}
\cvolyear{1993}
\cmonth{April}
\cyear{1993}
\cvolno{2}
\cpgs{253-287}
\title
New points of view in knot theory
\endtitle
\author
Joan S. Birman
\endauthor
\address
Department of Mathematics, 
Columbia University,
New York, New York 10027\endaddress
\ml jb\@math.columbia.edu\endml
\thanks
This research was supported in part by NSF Grant 
DMS-91-06584 and by the
 US-Israel Binational Science
Foundation Grant 89-00302-2\endthanks
\thanks This manuscript is an expanded version of an 
AMS-MAA Invited Address,
given on January 8, 1992, at the Joint Winter Meeting of 
the AMS and the MAA in
Baltimore, Maryland\endthanks
\date September 30, 1992\enddate
\endtopmatter

\document

\heading Introduction\endheading

In this article we shall give an account of certain 
developments  in knot
theory which followed upon the discovery of the Jones 
polynomial \cite {Jo3} in
1984. The focus of our account will be recent glimmerings 
of understanding of
the topological meaning of the new invariants. A second 
theme will be the
central role that braid theory has played  in the subject. 
A third will be the
unifying principles provided by representations of simple 
Lie algebras and
their universal  enveloping algebras. These choices in 
emphasis are our own.
They represent, at best, particular aspects of the 
far-reaching ramifications
that followed the discovery of the Jones polynomial.

We begin in \S 1 by discussing the topological 
underpinnings of that 
most famous of the classical knot invariants---the 
Alexander polynomial. 
It will serve as a model for the sort of thing one would 
like to see for the
Jones polynomial. Alexander's 1928 paper ends with a hint of
things to come, in the form of a crossing-change formula 
for his 
polynomial, and in \S 2 we discuss how related formulas 
made their
appearances in connection with the Jones polynomial and 
eventually
led to the discovery of other, more general knot and link 
polynomials. 
A more systematic description
of these ``generalized Jones invariants" is given in \S 3, 
via the 
theory of R-matrices. That is where braids enter the 
picture, because every generalized Jones invariant is 
obtained from a 
trace function on an ``R-matrix representation" of the 
family of
braid groups $\lbrace B_n; n = 1,2,3,\dots\rbrace.$  The 
mechanism for 
finding R-matrix representations of $B_n$ is the
theory of quantum groups. For this reason, 
the collection of knot and link invariants that they 
determine have been called 
{\it quantum group invariants \/}. We shall refer to them 
here
as either quantum group invariants or {\it generalized 
Jones invariants \/},
interchangeably. While the theory of R-matrices and their  
construction via
quantum groups gives a coherent and  uniform description 
of the entire class of
invariants, the underlying  ideas will be seen to be 
essentially
combinatorial in nature. Thus, by the end of \S 3 the 
reader should begin to
understand how it could happen that in 1990 topologists 
had a fairly coherent
framework for constructing  vast new families of knot and 
link invariants,
possibly even enough to classify unoriented knot and link 
types, without having
the slightest clue to the underlying topology.

In \S 4 we introduce an entirely new collection of 
invariants, which  arose out
of techniques pioneered by Arnold in singularity theory 
(see the introduction
to \cite {Arn1, Arn2}). The new invariants will be seen to 
have a solid basis
in a very interesting new topology, where one studies not 
a  {\it single \/}
knot, but a space of {\it all \/} knots. This  point of 
view was carried out
successfully for the case of knots  by Vassiliev \cite 
{V}. ({\it Remark\/}. 
For simplicity, we restrict our attention in this part of 
the review to knots.
The theory is,  at this writing, less well developed for 
links.)  Vassiliev's
knot invariants  are rational numbers. They lie in vector 
spaces $V_i$ of
dimension $d_i, i=1,2,3,\dots$, with invariants in $V_i$ 
having ``order" $i$.
On the other hand, quantum group invariants are Laurent 
polynomials over 
$\roman {Z}$, in a single variable $q$. The relationship 
between them is as
follows:  

\proclaim{Theorem 1} Let ${\Cal J}_q(\BK)$ be a quantum 
group
invariant for a knot type $\BK$. Let 
$${\Cal P}_x(\BK) = \sum_{i=0}^{\inf} u_i(\BK)x^i$$ 
be the power series over the rational numbers $Q$ obtained 
from ${\Cal J}_q(\BK)$ by setting $q=e^x$ and expanding 
the powers of $e^x$
in their Taylor series. Then the coefficient $u_i(\BK)$ of 
$x^i$ is $1$ if
$i=0$ and a Vassiliev invariant of order $i$ for each $i 
\geq 1$. 
\endproclaim

\np Thus Vassiliev's topology of the space of all knots 
suggests the topological 
underpinnings we seek for the quantum group knot invariants.

A crucial idea in the statement of Theorem 1 is that one 
must pass to the power
series representation of the knot polynomials before one 
can understand the
situation. This was first explained to the author and Lin 
by Bar Natan. Theorem
1 was first proved for special cases in  \cite {BL} and 
then generalized in
\cite {Li1}. In \S\S5--7 we describe a set of ideas  which 
will be seen to lead
to a new and  very simple proof of Theorem 1. First, in \S 
5, we review how
Vassiliev's  invariants,  like the Jones polynomial, the 
HOMFLY and Kauffman 
polynomials \cite {FHL, Kau1}, and the $\bold {G}_2$ 
polynomial \cite {Ku} can
be described by axioms and initial data.  (Actually, all 
of the quantum group
invariants admit such a description,  but the axioms can 
be very complicated;
so such a description would  probably not be 
enlightening). In \S 6 we
introduce braids into the Vassiliev setting, via a new 
type of object, the
monoid of singular braids. Remarkably, this monoid will be 
seen to map
homomorphically (we conjecture isomorphically) into the 
group algebra of the
braid group, implying that it is as fundamental an object 
in mathematics as the
braid group itself. This allows us to extend every 
R-matrix representation of
the braid group to the singular braid monoid. In \S7 we 
use R-matrices and
singular braids to prove Theorem 1. In \S8 we return to 
the topology,
discussing the beginning of a topological understanding of 
the quantum group
invariants.  We then discuss, briefly, a central problem: 
do we now know enough
to classify knots via their algebraic invariants? We will 
describe some of the
evidence which allows us to sharpen that question.

Our goal, throughout this review, is to present the material
in the most transparent and nontechnical manner possible 
in order to help
readers who work in other areas to learn as much as 
possible about
the state of the art in knot theory. Thus, when we give 
``proofs", they
will be, at best, sketches of proofs. We hope there will 
be enough detail
so that, say, a diligent graduate student who is motivated 
to
read a little beyond this paper will be able to fill in 
the gaps.  

Among the many topics which we decided deliberately to 
{\it exclude\/} from
this review for reasons of space, one stands out: it 
concerns the
generalizations of the quantum group invariants and 
Vassiliev invariants to
knots and links in arbitrary 3-manifolds, i.e., the 
program set forth by Witten
in \cite {Wi}. That very general program is inherently 
more difficult than the
special case of knots and links or simply of knots in the 
3-sphere. It is an
active area of research, with new discoveries made every 
day. We thought, at
first, to discuss, very briefly, the 3-manifold invariants 
of Reshetikhin and
Turaev \cite {RT} and the detailed working out of special 
cases of those
invariants by Kirby and Melvin \cite {KM}; however, we 
then realized that we 
could not include such a discussion and ignore  Jeffrey's 
formulas for the
Witten invariants of the lens spaces  \cite {Je}.  
Reluctantly, we
made the decision to restrict our attention to knots in 
3-space, but still, we
have given at best a restricted picture. For example, we 
could not do justice
to the topological constructions in  \cite {La} and in 
\cite {Koh1, Koh2} 
without making this review much longer than we wanted it 
to be, even though it
seems very likely that those constructions are closely 
related to the central
theme of this review.  We regret those and other omissions.

\heading 1. An introduction to knots and their Alexander 
invariants 
\endheading

We will regard a knot 
$K$ as the embedded image of an oriented circle 
$ S^1$ in oriented 3-space
$\Bbb{R}^3$ or $S^3$. 
If, instead, we begin with $\mu \geq 2 $ copies of $S^1$, 
the image (still called $K$), is a {\it link \/}. 
Its {\it type \/} $\BK$ is the 
topological type of the pair $(S^3,K)$, 
under homeomorphisms which preserve 
orientations on both $S^3$ and $K$. Knot types do not 
change if we replace $S^3$ by
$\Bbb{R}^3,$ because every homeomorphism of $S^3$ is 
isotopic to one which fixes one point, 
and that point may be chosen to be the point at infinity.

We may visualize a knot via a  {\it diagram \/}, i.e., a 
projection of $K
\subset \Bbb{R}^3$ onto a generic  $\Bbb{R}^2 \subset 
\Bbb{R}^3$,  where the
image is decorated to distinguish overpasses from 
underpasses, for example, as
for highways on a map. Examples are given in Figure 1.  
The example in Figure
1(a) is  the unknot, represented by a round planar circle. 
Figure 1(b) shows a 
{\it layered \/} diagram of a knot, i.e., one which has 
been drawn, using an
arbitrary but fixed starting point (which in Figure 1(c) 
is the tip of the
arrow) without the use of an eraser, so that the first 
passage across a double
point in the projection is always an overpass. We leave it 
to the reader to
prove that a layered diagram with $\mu$ components always 
represents $\mu$
unknotted,  unlinked circles. From this simple fact it 
follows that {\it any
\/} diagram of {\it any \/} link may be systematically 
changed to a diagram for
the unlink on the same number of components by finitely 
many crossing changes. 

\midspace{23pc}
\caption{{\smc Figure 1}. {\rm Diagrams of the unknot.}}

The diagram in Figure 1(c) was chosen to illustrate the 
subtleties of knot
diagrams. It too represents the unknot, but it is not 
layered. This example was
constructed by Goeritz  \cite {Go} in the mid-1930s. At 
that time it was known
that a finite number of repetitions of Reidemeister's 
three famous moves,
depicted (up to obvious symmetries and variations) in 
Figure 2, suffice to 
take any one diagram of a knot to any other. Notice that 
Reidemeister's moves
are ``local" in the sense that they are restricted to 
regions of the diagram
which contain at most three crossings.  If one  defines 
the {\it complexity \/}
of a knot diagram to be its crossing number, a natural 
question is whether
there is a set of moves that preserve or reduce complexity 
and that, applied
repeatedly, suffice to reduce  an arbitrary diagram of an 
arbitrary knot to one
of minimum  crossing number. The diagram in Figure 1(c) 
effectively ended that
approach to the subject, because the eight moves which had 
by then been
proposed as augmented Reidemeister moves did not suffice 
to simplify this
diagram. 

Some insight may be obtained into the question by 
inspecting the ``handle move"
of Figure 3. Note that the crossing number in the diagram 
of Figure 1(c) can be
reduced by an appropriate handle move. The region that is 
labeled $X$ in Figure
3 is arbitrary. The handle move clearly decreases crossing 
number, but a few
moments\vadjust{\fighere{12.5pc}\caption{{\smc Figure 2}.
{\rm Reidemeister moves.}}} of thought should convince
the reader that if one tries to factorize it
into a product of Reidemeister moves, for any sufficiently 
complicated choice
of $X$, it will be necessary to use the second 
Reidemeister move repeatedly in
order to create regions in which the third move may be 
applied. The same sort
of reasoning makes it unlikely that any set of local moves 
suffices to reduce
complexity. Indeed, if we knew such a set of moves, we 
would have the beginning
of an algorithm for solving the knot problem, because 
there are only finitely
many knot diagrams with fixed crossing number.

\topspace{6pc}
\caption{{\smc Figure 3}. {\rm A handle move.}} 

In Figure 4 we have given additional examples, taken from 
the beginning  of a
table compiled at the end of the nineteenth century by the 
physicist  Peter
Guthrie Tait and coworkers, as part of a systematic effort 
to classify knots.
The knots in the table are listed in order of their 
crossing number, so that,
for example, $7_6$ is the 6th knot that was discovered 
with 7 crossings. The
part of the table which we have shown includes all prime 
knots with at most $7$
crossings, up to the symmetries defined by reversing the 
orientation on either
$K$ or $S^3$. (We have shown the two trefoils, for reasons 
which will become
clear shortly.) The tables, which eventually included all 
prime knots defined
by diagrams with at most 13 crossings, were an ambitious 
undertaking. (Aside:
Yes, it is true. Knots, like integers, have a 
decomposition into an
appropriately defined product of prime knots, and this 
decomposition is unique
up to order.) Their clearly stated goal \cite {Ta} was to  
uncover the
underlying principles of knotting, but to the great  
disappointment of all
concerned they did not even succeed in revealing a single 
knot type invariant
which could be computed from a diagram. Their importance, 
to this day, is due
to the fact that they provide a rich source of examples 
and convincing evidence
of both the beauty and subtlety of the subject.

\topspace{37.5pc}
\caption{{\smc Figure} 4.}

Beneath each of the examples in Figure 4 we show two 
famous  invariants,
namely, the Alexander polynomial $A_q(\BK)$ and the 
one-variable Jones 
polynomial $J_q(\BK)$. Both are to be regarded as Laurent 
polynomials in $q$,
the series of numbers representing the sequence of 
coefficients, the bracketed
one being the constant term. Thus the lower sequence 
(1,-1,1,-2, \lbrack
2\rbrack , -1,1)  beneath the knot $6_1$ shows that its 
Jones polynomial is
$q^{-4} - q^{-3} + q^{-2} - 2q^{-1} + 2 - q + q^2$. The 
richness of structure
of both invariants is immediately clear from the sixteen 
examples in our table.
There are no duplications (except for the Alexander 
polynomials of the two
trefoils, which we put in deliberately to make a point).  
The arrays of
subscripts and superscripts, as well as the roots (which 
are not shown) and the
poles of the Jones invariant (about which almost nothing 
is known), suggest
that the polynomials could encode interesting properties 
of knots. Notice that
$A_q(\BK)$ is symmetric, i.e., $A_q(\BK) = \bold 
{A}_{q^{-1}}(\BK)$, or,
equivalently, the array of coefficients is palendromic. On 
the other hand,
$J_q(\BK)$ is not. Both polynomials take the value 1 at 
$q=1$.  ({\it
Remark\/}. 
The Alexander polynomial is actually only determined up to
$\pm$multiplicative powers of $q$, and we have chosen to 
normalize it to stress
the symmetry and so that its value at $1$ is $+1$ rather 
than $-1$.)

The knots in our table are all invertible; i.e., there is 
an isotopy of 3-space
which takes the oriented knot to itself with reversed 
orientation. The first
knot in the tables that fails to have that property is 
$8_{17}$. Neither the
Alexander nor the Jones polynomials changes when the 
orientation on the knot is
changed.  We shall have more to say about noninvertible 
knots in \S 8.

While we understand the underlying  topological meaning of 
$A_q(\BK)$ very
well, we are only beginning to  understand the topological 
underpinnings of
$J_q(\BK)$.  To begin to explain the first assertion, let 
us go back to one of
the earliest problems in knot theory: to what extent does 
the topological type
$\bold{X}$ of the complementary space $X = S^3 - K$ and/or 
the isomorphism
class $\bold{G}$ of  its fundamental group  $G(K) = 
\pi_1(X,x_0)$ suffice to
classify knots?  The trefoil knot is almost everybody's 
candidate for the
simplest example of a nontrivial knot,  so it seems 
remarkable that, not long
after the discovery of the fundamental group of a 
topological space, Max Dehn
\cite {De} succeeded in proving that the trefoil knot and 
its mirror image had 
isomorphic groups, but that their knot types were 
distinct. Dehn's proof is
very ingenious! This was the beginning of a long story, 
with many contributions
(e.g., see \cite {Sei, Si, CGL}) which reduced repeatedly 
the number of
distinct knot types which could have homeomorphic 
complements and/or isomorphic
groups, until it was finally proved, very recently, that 
(i) $\bold {X}$
determines $\BK$ (see \cite {GL}) and (ii) if $\BK$ is 
prime, then $\bold{G}$
determines $\BK$ up to unoriented equivalence \cite {Wh}. 
Thus there are at
most four distinct oriented  prime knot types which have 
the  same knot group
\cite {Wh}. This fact will be important to us shortly.

The knot group $\bold{G}$ is finitely presented; however, 
it is infinite,
torsion-free, and (if $\BK$ is not the unknot) nonabelian. 
Its isomorphism
class is in general not easily understood via a direct 
attack on the problem.
In such circumstances, the obvious thing to do is to pass 
to the abelianized
group, but unfortunately $G/[G,G]  \cong H_1(X;\roman 
{Z})$ is infinite cyclic
for all knots, so it is of no use in distinguishing knots. 
Passing to the
covering space $\widetilde X$  which belongs to $[G,G]$, 
we note that there is a
natural action of the cyclic group $G/[G,G]$  on 
$\widetilde X$ via covering
translations. The action makes the homology group 
$H_1(\widetilde X;\roman
{Z})$ into a $\roman {Z}[q, q^{-1}]$-module,  where $q$ is 
the generator of
$G/[G,G]$. This module turns  out to be finitely 
generated. It's the famous
{\it Alexander module \/}. While the ring $\roman {Z}[q, 
q^{-1}]$ is not a PID,
relevant aspects of the theory of modules over a PID apply 
to $H_1(\widetilde
X;\roman {Z})$. In particular, it splits as a direct sum 
of cyclic modules, the
first nontrivial one being $\roman 
{Z}[q,q^{-1}]/A_q(\BK)$. Thus $A_q(\BK)$ is
the generator of the ``order ideal", and the smallest 
nontrivial torsion
coefficient in the module $H_1(\widetilde X)$. In 
particular,  $A_q(\BK)$ is
very clearly an invariant of the knot group.

We regard the above description of $A_q(\BK)$ as an 
excellent model  for what
we might wish in a topological invariant of knots. We know 
precisely what it
detects, and so we also know precisely what it {\it fails 
\/} to detect. For
example, it turns out that $\pi _1(\widetilde X)$ is 
finitely generated if and
only if $X$ has the structure of a surface  bundle over 
$S^1$, but there is no
way to tell definitively from $A_q(\BK)$ whether $\pi 
_1(\widetilde X)$ is or
is not finitely generated. On the other hand, if a surface 
bundle structure 
exists, the genus of the surface is determined by 
$A_q(\BK)$.  The polynomial
$A_q(\BK)$ also generalizes in many ways. For example, 
there are Alexander
invariants of links, also additional Alexander invariants 
when the Alexander
module has more than one torsion coefficient. Moreover, 
the entire theory
generalizes naturally to higher dimensional knots, the 
generalized invariants 
playing a central role in that subject.

Returning to the table in Figure 4, we remark that when a 
knot  is replaced by
its mirror image (i.e., the orientation on $S^3$ is 
reversed), the Alexander and
Jones polynomials $A_q(\BK)$ and $J_q(\BK)$ go over to 
$A_{q^{-1}}(\BK)$ and 
$J_{q^{-1}}(\BK)$ respectively. As noted earlier, 
$A_q(\BK)$ is invariant under
such a change, but from the simplest possible example, the 
trefoil knot, we see
that $J_q(\BK)$  is not. Now recall that $\bold {G}$ does 
{\it not \/} change
under changes in the  orientation of $S^3$. This simple 
argument shows that
$J_q(\BK)$ cannot be a group  invariant!  Since there are 
at most four distinct
knot types that share the same knot group $\bold{G}$,  a 
first wild guess would
be that $J_q(\BK)$, which does detect changes in the 
ambient space orientation
(but not in knot orientation), classifies unoriented knot 
types;   but this 
cannot be true because \cite {Kan} constructs examples of 
infinitely many
distinct prime knot types with the same Jones polynomial. 
Thus it seems
interesting indeed to ask about the underlying topology 
behind the Jones
polynomial. If it is not a knot group invariant, what can 
it be?  We will begin
to approach that problem by a circuitous route, taking as 
a hint the central
and very surprising role of ``crossing-change formulas"  
in the subject. 

\heading 2. Crossing changes
\endheading

A reader who is interested in the history of mathematics 
will find,  on
browsing through several of Alexander's Collected Works, 
that many of his 
papers end with an intriguing or puzzling comment or 
remark which, as it turned
out with the wisdom of hindsight, hinted at future 
developments of the
subject.  For example, in his famous paper on braids \cite 
{Al1}, which we will
discuss in detail in \S 6, he proves that  every knot or 
link may be
represented as a closed braid. He then remarks (at the end 
of the paper) that
this yields a construction for describing 3-manifolds via 
their fibered knots;
however, he did it long before anyone had considered the 
concept of a fibered
knot!  Another example that is of direct interest to us 
now occurs in \cite
{Al2}, where he reports on the discovery of the Alexander 
polynomial. In
equation (12.2) of that paper  we find observations on the 
relationship between
the Alexander polynomials of three  links: $K_{p_+}, 
K_{p_-},$ and $K_{p_0}$, 
which are defined by diagrams  that are identical outside 
a neighborhood of a
particular double point $p$,  where they differ in the 
manner indicated in
Figure 5. 

The formula which Alexander gives is:
$$ A_q(\BK_{p_+}) - A_q(\BK_{p_-}) = (\sqrt q - 1/{\sqrt 
q}) A_q(\BK_{p_0}).
\tag 1a $$
This formula passed unnoticed for forty years. (We first 
learned about
Alexander's version of it in 1970 from Mark Kidwell.) 
Then, in 1968 it was
rediscovered, independently, by John Conway \cite {C}, who 
added a new
observation: If you require, in addition to (1a), that:
$$ A_q(\bold {O}) = 1, \tag 1b $$
where ${\bold {O}}$ is the unknot, then (1a) and (1b) 
determine $A_q(\BK)$ on
all knots, by a 
recursive\vadjust{\fighere{5.5pc}\caption{{\smc Figure 5}. 
{\rm Related link diagrams.}}}
procedure. To see this, the first thing to observe 
is that if $\bold{O}_{\mu}$ is the $\mu$-component unlink, 
then we may find related diagrams in which 
$K_{p_+}, K_{p_-},$ and $K_{p_0}$ represent 
$\bold {O}_{\mu}, \bold {O}_{\mu}$, and $\bold {O}_{\mu + 
1}$ respectively, 
as in Figure 6.

\topspace{4pc}
\caption{{\smc Figure 6}. {\rm Three related diagrams for 
the unlink.}}

This fact, in conjunction with (1a), implies that 
$A_q(\bold{O}) = 0$ if
$\mu \geq 2.$ Next, recall (Figure 1b) that any diagram 
$K$ for any knot type
$\BK$ may be changed to a layered diagram which represents 
the unknot or unlink
on the same number of components, by appropriate crossing 
changes. Induction on
the number of crossing changes to the unlink then 
completes the proof of 
Conway's result. This led him to far-reaching 
investigations of the
combinatorics of knot diagrams.

While the Jones polynomial was discovered via braid theory 
(and we shall have
more to say about that shortly), Jones noticed, very early 
in the game, that
his polynomial also satisfied a crossing-change formula, 
vis:
$$ q^{-1}J_q(\BK_{p_+}) - qJ_q(\BK_{p_-}) = 
(\sqrt q - 1/{\sqrt q}) J_q(\BK_{p_0}). \tag  2a $$
which, via Figure 6 and layered diagrams, may be used, in 
conjunction
with the initial data 
$$ J_q(\bold O) = 1, \tag  2b $$
to compute $J_q(\BK)$ for all knots and links. Motivated 
by the
similarity between (1a)--(1b) and (2a)--(2b), several 
authors \cite {LM,
Ho, PT} were led to consider a more general 
crossing-change formula, which for our purposes may be 
described as an infinite
sequence of crossing-change formulas:
$$ q^{-n}H_{q,n}(\BK_{p_+}) - q^nH_{q,n}(\BK_{p_-}) = 
(\sqrt q - 1/{\sqrt q}) H_{q,n}(\BK_{p_0}), \tag  3a $$
$$ H_{q,n}(\bold O) = 1, \tag 3b $$
where $n \in \roman{Z}$. It turns out that (3a)--(3b) 
determine an  infinite
sequence of one-variable polynomials which in turn extend 
uniquely to give a
two-variable invariant which has since become known as the 
HOMFLY polynomial
\cite {FHL}. Later, (1a)--(3b) were replaced by a more 
complicated family of
crossing-change formulas, yielding the Kauffman polynomial 
invariant of knots
and links \cite {Kau1}. A unifying principle  was 
discovered (see \cite {Re2})
which yielded still further invariants, for example, the 
$G_2$ invariant of
embedded knots, links, and graphs \cite {Ku}. Later, 
crossing-change formulas
were used to  determine other polynomial invariants of 
knotted graphs \cite
{Y1} as well.   

At this time we know many other polynomial invariants of  
knots, links, and
graphs. In principle, all of them can be defined  via 
generalized
crossing-changes, together with initial data. In general, 
a particular
polynomial will be defined by a family of equations which 
are like (2a) and
(3a). It is to be expected that each such equation will 
relate the invariants
of knots which are defined by diagrams which differ in 
some specified way, in a
region which has a fixed number of incoming and outgoing 
arcs. The ways in
which they differ will be more complicated than a simple 
change in a crossing
and ``surgery"  of the crossing. Thus we have a conundrum: 
on the one hand,
knot and\ link diagrams and crossing-change formulas 
clearly have much to do
with the  subject; on the other hand, their role is in 
many ways puzzling,
because we do not seem to be learning as much as we might 
expect to learn about
diagram-related invariants from the polynomials. 

We make this last assertion explicit. First, let us define 
three diagram-related 
knot invariants: 

(i) the minimum crossing number $c(\BK)$ of a knot,

(ii) the minimum number of crossing changes $u(\BK)$, 
i.e., to the unknot, and

(iii) the minimum number of Seifert circles $s(\BK)$,  

\np where each of these invariants is to be minimized over 
all possible knot
diagrams. ({\it Aside\/}. For the reader who is unfamiliar 
with the concept of
Seifert circles, we note that by \cite {Y2} the integer 
$s(\BK)$ may also be
defined to be the braid index of a knot or link, i.e., the 
smallest integer $s$
such that $\BK$ may be represented as a closed $s$-braid. 
If the reader is also
not familiar with closed braids, he or she might wish to 
peek ahead to \S3.)
All of them satisfy inequalities which are  detected by 
knot polynomials. For
example, the Morton-Franks-Williams inequality places 
upper and lower  bounds
on $s(\BK)$ \cite {Mo2, FW}. Also, the Bennequin 
inequality, recently
proved by Menasco \cite {Me}, gives a lower bound for 
$u(\BK)$ which can be
detected by the one-variable Jones polynomial. As another 
example, the
one-variable Jones polynomial was a major tool in the 
proof of the Tait
conjecture (see \cite {Kau2, Mu}), which relates to the 
definitive
determination of $c(\BK)$ for the special case of 
alternating knots which are
defined by alternating diagrams. On the other hand, there 
are examples which
show that none of the inequalities mentioned above can be 
equalities for all
knots.  Indeed, at this writing we do not have a 
definitive method for
computing any of these very intuitive diagram-related knot 
and link-type
invariants.

\heading 3. R-Matrix representations of the braid group
\endheading

One of the very striking successes of the past eight years 
is that, after a
period during which new polynomial invariants of knots and 
links were being
discovered at an alarming rate (e.g., see \cite {WAD}), 
order came out of chaos
and a unifying principle emerged which gave a fairly 
complete  description of
the new invariants; we study it in this section. We remark 
 that this
description may be given in at least two mutually 
equivalent ways:  the first
is via the algebras of \cite {Jo4, BW, Kal} and the 
cabling construction of 
\cite {Mu, We2}; the second is  via the theory of 
R-matrices, as we shall do
here.

\topspace{16.5pc}
\caption{{\smc Figure} 7. {\rm Braids.}}

Our story begins with the by-now familiar notion of an 
$n$-braid \cite {Art}.
See Figure 7(a) for a picture, when $n=3$. 
Our $n$-braid is to be regarded as living in a slab of 
3-space $ \Bbb{R}^2
\times  I \subset  \Bbb{R}^3.$  It consists of $n$ 
interwoven oriented strings
which join $n$ points, labeled $1,2,\dots,n$, in the plane 
$ \Bbb{R}^2 \times
\lbrace 0 \rbrace$ with corresponding points in $ 
\Bbb{R}^2 \times \lbrace 1
\rbrace$,  intersecting each intermediate plane $ 
\Bbb{R}^2 \times \lbrace t
\rbrace$ in exactly $n$ points. Two braids are equivalent 
if there is an
isotopy of one to the other which preserves the initial 
and end point of each
string, fixes each plane $ \Bbb{R}^2 \times \lbrace t 
\rbrace$ setwise, and
never allows two strings to intersect. Multiplication is 
by  juxtaposition,
erasure of the middle plane, and rescaling. Closed braids 
are obtained from
(open) braids by joining the initial point of each strand, 
in $ \Bbb{R}^3$, to
the corresponding end point, in the manner  illustrated in 
Figure 7(b), so that
if one thinks of the closed braid as wrapping around the 
$z$-axis, it meets
each plane $\theta=$constant in exactly $n$ points. A 
famous theorem of Alexander
\cite {Al1} asserts that every knot or link may be  so 
represented, for some
$n$. ({\it Remark\/}. Lemma 1 of \S6 sketches a  
generalization of Alexander's
original proof.) In fact, if $K$ is our oriented link,  we 
may choose any
oriented $ \Bbb{R}^1 \subset  \Bbb{R}^3$ which is disjoint 
 from $K$ and modify
$K$ by isotopy so that this  copy of $\Bbb{R}^1$ is the 
braid axis.

The braid group $B_n$ is generated by the elementary braids
$\sigma_1,\dots,\sigma_{n-1}$ illustrated in Figure 8. For
example, the braid in Figure 7(a) may be described, using
the generators given in Figure 8, by the {\it word}\ \ 
$\sigma_1^{-2} \sigma_2 \sigma_1^{-1}$. 
Defining relations in $B_n$ are:
$$ \sigma_i \sigma_j = \sigma_j \sigma_i \quad \text{if }
\mid i-j \mid \geq 2,  \tag 4a$$

$$ \sigma_i \sigma_j \sigma_i = \sigma_j \sigma_i \sigma_j
\quad \text{if } \mid i-j \mid = 1,  \tag 4b$$
\np We will refer to these as the {\it braid relations}.

Let $B_{\infty}$ be the disjoint union of the braid groups 
$ B_1,  B_2,\dots.$ 
Define $\beta,  \beta^{\ast} \in B_{\infty}$ to be {\it
Markov-equivalent \/} if the closed braids $\hat {\beta}, 
\hat {\beta}^{\ast}$
which they define represent the same link type. 
Markov's theorem, announced\vadjust{\fighere{9.5pc}
\caption{{\smc Figure 8}. {\rm Elementary braids, singular 
braids, and
tangles.}}}
in \cite {Mar} and proved forty years later in \cite {Bi}, 
asserts that Markov
equivalence is generated by conjugacy in each $B_n$ and 
the map $B_n
\rightarrow B_{n+1}$ which takes a word 
$W(\sigma_1,\dots\sigma_{n-1})$ to
$W(\sigma_1,\dots\sigma_{n-1}) \sigma_n^{\pm 1}$. We call 
the latter {\it
Markov's second move \/}.  The examples in Figures 7(a), 
7(b) may be used to
illustrate Markov equivalence. The 3-braid $\sigma_1^{-2} 
\sigma_2
\sigma_1^{-1}$ shown there is conjugate in $B_3$ to 
$\sigma_1^{-3} \sigma_2$,
which may be modified to the 2-braid $\sigma_1^{-3}$ using 
Markov's second
move, so that $\sigma_1^{-3} \in  B_2$ is Markov 
equivalent to $\sigma_1^{-2}
\sigma_2 \sigma_1^{-1} \in B_3$. 

Let $\{ \rho_n:B_n \rightarrow \GL _{m_n}(\Cal{E}); 
n=1,2,3,\dots\}$ be a
family of matrix representations of $B_n$ over a ring 
$\Cal{E}$ with 1. Let
$i_n:B_n \rightarrow B_{n+1}$ be the natural inclusion map 
which takes $B_n$ to
the subgroup of $B_{n+1}$ of braids on the first $n$ 
strings. A linear function
$\TR:\rho_n (B_n) \rightarrow \Cal{E}$ is said to be a 
{\it Markov trace \/}
if: 

(i) $\TR(1)=1$,

(ii) $\TR(\rho_n (\alpha \beta)) = \TR(\rho_n (\beta 
\alpha))\ \ \
\forall \alpha, \beta \in B_n$,

(iii) $\exists z \in \Cal{E}$ such that if $\beta \in 
i_n(B_n)$, then
$\TR(\rho_{n+1}(\beta \sigma_n^{\pm 1})) = 
(z)(\TR(\rho_{n}(\beta))$.

\

\np In particular, by setting $\beta = 1$ in (iii) we see 
that $z =
\TR(\rho_n(\sigma_{n-1}^{\pm 1}))$. In view of Markov's 
theorem, it is
immediate that every family of representations of $B_n$ 
which supports a Markov
trace determines a link type invariant  
$F(\bold{K}_{\beta})$ of the link type
$\bold{K}_{\beta}$ of the  closed braid $\hat \beta$, 
defined by 
$$F(\bold{K}_{\beta}) = z^{1-n}\TR(\rho_n(\beta)). \leqno 
(5)$$

Later, we will need to ask what happens to (5) if we 
rescale the
representation, so we set the stage for the modifications 
now. Suppose that
$\rho = \rho_n$ is a representation of $B_n$ in $\GL 
_m(\Cal{E})$ which
supports a function  $\TR$ which satisfies (i)--(iii). Let 
$\gamma$ be any
invertible element in $\Cal{E}$. Then we may  define a new 
representation 
$\rho^{\prime}$ by the rule $\rho^{\prime}(\sigma_i) = 
\gamma  \rho(\sigma_i)$
for each $i = 1,\dots,n$. It is easy to see that 
$\rho^{\prime}$ is a
representation  if and only if $\rho$ is, because (4a) and 
(4b) are satisfied
for one if and  only if they are satisfied for the other. 
Properties (i) and
(ii) will continue to hold if we replace $\rho$ by 
$\rho^{\prime}$, but (iii)
will be modified because the traces of $\sigma_i$ and 
$\sigma_i^{-1}$ will not
be the same in the new representation. To determine the 
effect of the change on
equation (5), let $z^{\prime} = 
\TR(\rho^{\prime}(\sigma_i)) = \gamma \
\TR(\rho(\sigma_i)) = \gamma z$, so $z = \gamma ^{-1} 
z^{\prime}$. Choose any
$\beta \in B_n$. Then $\beta$ may be expressed as a word
$\sigma_{\mu_1}^{\varepsilon_1} 
\sigma_{\mu_2}^{\varepsilon_2} \cdots
\sigma_{\mu_r}^{\varepsilon_r}$  in the generators, where 
each $\varepsilon_j =
\pm 1$. Let $\varepsilon (\beta) = \sum_i^r 
\varepsilon_i$. Then equation (5)
will be replaced by
$$F(\bold{K}_{\beta}) =
\gamma^{-\varepsilon(\beta)}(z^{\prime})^{1-n}\TR(\rho_n^{%
\prime}(\beta)). 
\leqno (5^{\prime})$$ 
Thus the existence of a link type invariant will not be 
changed
if we rescale the representation by introducing the 
invertible element
$\gamma$. Similar considerations apply if we rescale the 
trace function. For
example, instead of requiring that $\TR(1) = 1$, we could 
require that if $1_n$
denotes the identity element in $B_n$, then $\TR(1_n) = 
z^{n-1}$, in which case
the factor $z^{1-n}$ in equation (5) would vanish.
  
The invariant which is defined by equation (5) or 
(5$^{\prime}$) depends
directly on the traces of matrices in a finite-dimensional 
matrix
representation $\rho_n$ of  the braid group $B_n$. Any 
such representation is
determined by its values on the generators $\sigma _i$, 
and since these are
finite-dimensional matrices they satisfy their 
characteristic polynomials,
yielding polynomial identities. From the point of view of 
this section, such
identities are one source of the  crossing-change formulas 
we mentioned
earlier, in \S2. Another source will be trace identities, 
which always exist in
matrix groups. (See, e.g., \cite {Pr}.) In general one 
needs many such
identities (i.e., polynomial equations satisfied by  the 
images of various
special braid words) to obtain axioms which suffice to 
determine a link type
invariant.

In the summer of 1984, almost eight years ago to the day 
from this writing, this
author met with Vaughan Jones to discuss possible 
ramifications of a discovery
Jones had made in \cite {Jo2} of a new family of matrix 
representations of
$B_n$, in conjunction with his earlier studies of type 
II$_1$ factors and their
subfactors in von Neumann algebras. Before that time, 
there was essentially
only one matrix representation of the braid group which 
had a chance of being
faithful and had been studied in any detail---the Burau 
representation (see
\cite {Bi}). The knot invariant which was associated to 
that representation was
the Alexander polynomial. Jones had shown that his 
representations contained
the Burau representation as a proper summand. Thus his 
representations were new
and interesting. Also, they supported an interesting trace 
function.  As  it
developed \cite {Jo3}, his trace functions were Markov 
traces. If $\rho _n$ is 
taken to be the Jones's representation of $B_n$, the 
invariant we called
$F(\BK_{\beta})$ in (5) is the one-variable Jones 
polynomial $J_q(\BK)$. 
  
\subheading{Wandering along a bypath}
The reader who is pressed for time may wish to omit this 
detour.  We interrupt
our main argument to discuss how the Jones representations 
were discovered,
because it is very interesting to see how an almost chance 
discovery of an
unexpected relationship between two widely separated areas 
of mathematics had
ramifications which promise to keep mathematicians busy at 
work for years to
come! Let $M$ denote a von Neumann algebra. Thus $M$ is an 
algebra of  bounded
operators  acting on a Hilbert space $\Cal{H}$. The 
algebra $M$ is called a
factor if its  center  consists only of scalar multiples 
of the identity. The
factor is type II$_1$ if  it admits a linear functional 
$\TR:M \rightarrow
\Bbb{C}$, which satisfies:

(i) $\TR(1) = 1$,

(ii) $\TR(xy) = \TR(yx) \ \ \ \forall x,y \in M$, 

\np and a positivity condition which shall not concern us 
here. It is known
that the trace is unique, in the sense that it is   the 
only linear form
satisfying (i) and (ii). An old discovery of  Murray and 
von Neumann was that
factors of type II$_1$ provide a type of ``scale"  by 
which one can measure the
dimension $\roman{dim}_M\Cal{H}$ of the Hilbert space 
$\Cal{H}$.  The notion of
dimension  which occurs here generalizes the familiar 
notion of integer-valued
dimensions,  because for appropriate $M$ and $\Cal{H}$ it 
can be any
nonnegative real  number or infinity.	The starting point 
of Jones's work was
the following question:  if $M_1$ is a  type II$_1$ factor 
and if  $M_0 \subset
M_1$ is a subfactor,  is there any restriction on the  
real numbers $\lambda$
which occur as the ratio $\lambda =
\roman{dim}_{M_0}\Cal{H}/\roman{dim}_{M_1}\Cal{H}$?
   
The question has the flavor of questions one studies in 
Galois theory. On the 
face of it, there was no reason to think that $\lambda$ 
could not take on any
value in  [1,$\infty$], so Jones's answer came as a 
complete surprise. He called
$\lambda$ the  {\it index \/}  [$M_1:M_0$] of $M_0$ in 
$M_1$ and proved a
type of rigidity theorem about it:

\proclaim {The Jones Index Theorem}  If $M_1$ is a 
$\roman{II}_1$ factor and 
$M_0$ a subfactor, then the  possible values of the index 
$\lambda$ are
restricted to $[4,\infty] \cup [4\cos^2(\pi /p)]$,  where 
$p \geq 3$ is a
natural number. Moreover, each such real number occurs for 
some pair $M_0,
M_1$.\endproclaim

We now sketch the idea of Jones's proof, which is to be 
found in \cite {Jo1}.
Jones begins with the type II$_1$ factor $M_1$ and a 
subfactor  $M_0$. There is
also a tiny bit of additional structure:  In this setting 
there exists a map
$e_1:M_1 \rightarrow M_0,$  known as the {\it conditional 
expectation \/} of
$M_1$ on  $M_0$. The map $e_1$ is a projection, i.e., 
$e_1^2 = e_1$.  His first
step is to prove that the ratio $\lambda$ is independent 
of the choice of  the
Hilbert space $\Cal{H}$. This allows him to choose an 
appropriate $\Cal{H}$  so
that the  algebra $M_2 = \langle M_1,e_1\rangle$ generated 
by $M_1$ and $e_1$ 
makes sense. He then investigates $M_2$ and proves that it 
is another  type
II$_1$ factor, which contains $M_1$ as a  subfactor; 
moreover, $|M_2:M_1| =
|M_1:M_0| = \lambda$.  Having in hand another type II$_1$ 
factor, i.e., $M_2$
and its subfactor $M_1$,  there is also a trace on $M_2$ 
which (by the
uniqueness of the trace)  coincides with the trace on 
$M_1$ when  it is
restricted to $M_1$ and another conditional expectation  
$e_2:M_2 \rightarrow
M_1$. This  allows Jones to iterate the construction and 
to build algebras 
$M_1,M_2, \dots$ and from them a family of algebras:
$$ \Cal{A}_n = \langle 1,e_1,...,e_{n-1} \rangle \subset 
M_n,\qquad 
n=1,2,3,\dots\.$$
We now replace the $e_k$\<'s by a new set of generators 
which 
are units, defining
$g_k = qe_k - 1 + e_k$, where 
$(1-q)(1-q^{-1}) = 1/\lambda$.
The $g_k$\<'s generate $\Cal{A}_n$, because the $e_k$\<'s 
do,  
and we can solve for the $e_k$\<'s in terms of the 
$g_k$\<'s.  Thus 
$$\Cal{A}_n = \Cal{A}_n(q) =  \langle 
1,g_1,\dots,g_{n-1}\rangle, $$
and we have a tower of algebras, ordered by inclusion:
$$\Cal{A}_1(q) \subset \Cal{A}_2(q) \subset \Cal{A}_3(q) 
\subset \cdots.$$ 
The parameter $q$, which replaces the index $\lambda$, is 
the quantity now 
under investigation. 	 
The $g_i$\<'s turn out to be 
invertible and to satisfy the braid relations (4a)--(4b), so
that there is a homomorphism from $B_n$ to $\Cal{A} _n$, 
defined by mapping
the elementary braid $\sigma_i$ to $g_i$. The parameter 
$q$ is woven into the 
construction of the tower. Defining relations in 
$\Cal{A}_n(q)$ 
depend upon $q$, for example, the relation $g_i^2 = 
(q-1)g_i + q$ holds.
Recall that, since $M_n$ is type II$_1$, it supports a 
unique trace, 
and, since $\Cal{A}_n$ is a subalgebra,
it does too, by restriction. This trace is a 
Markov trace! Jones's proof of the Index Theorem is 
concluded when he shows that the
infinite sequence of algebras  $\Cal{A}_n(q)$, with the 
given trace, could not exist
if $q$ did not satisfy the stated restrictions. 

Thus the ``independent variable" in the Jones polynomial 
is essentially the
index of a type II$_1$ subfactor in a type II$_1$ factor! 
Its discovery opened
a  new chapter in knot and link theory. 

\subheading{Back to the main road} We now describe a 
method, discovered by
Jones (see  the discussion of vertex models in \cite 
{Jo5})  but first worked
out in full detail by Turaev in  \cite {Tu}, which can be 
applied  to give, in
a unified setting, every generalized Jones invariant via a 
Markov trace on an
appropriate matrix representation of $B_n$. As before, 
$\Cal{E}$ is a ring with
1.  Let $V$ be a free $\Cal{E}$-module of rank $m \geq 1$. 
For each $n \geq 1$
let $V^{\otimes n}$ denote the $n$-fold tensor product $V 
\otimes_{\Cal{E}} \dots
\otimes_{\Cal{E}}V$. Choose a basis $v_1,\dots,v_m$ for 
$V$, and choose a
corresponding basis $\lbrace v_{i_1} \otimes \dots \otimes 
v_{i_n} ; 1
\leq i_1,\dots ,i_n \leq m \rbrace$ for $V^{\otimes n}$. 
An $\Cal{E}$-linear
isomorphism $f$ of $V^{\otimes n}$ may then be represented 
by an 
$m^n$-dimensional matrix $(f_{i_1 \cdots i_n}^{j_1 \cdots 
j_n})$ over
$\Cal{E}$, where the $i_k$\<'s (resp. $j_k$\<'s) are row 
(resp.\ column)
indices.

The family of representations of $B_n$ which we now 
describe have a very special form.
They are completely determined by an $\Cal{E}$-linear 
isomorphism $\roman{R}:V^{\otimes 2}
\rightarrow V^{\otimes 2}$ with matrix $\lbrack 
\roman{R}_{i_1i_2}^{j_1j_1} \rbrack$ as above.
Let $I_V$ denote the identity map on the vector space $V$. 
The representation 
$\rho _{n,\roman{R}}:B_n \rightarrow \GL _{m^n}(\Cal{E})$ 
that we need is defined by 
$$\rho _{n,\roman{R}}(\sigma _i) = I_V \otimes \dots 
\otimes I_V 
\otimes \roman{R}
\otimes I_V \otimes \dots \otimes I_V \leqno (6)$$  
where $\roman{R}$ acts on the $i$\<th and $(i+1)$\<st 
copies of $V$ in $V^{\otimes n}$. Thus, if we know how R 
acts on $V^{\otimes 2}$, we
know $\rho _{n,\roman{R}}$ for every natural number $n$. 

What properties must R satisfy for $\rho _{n,\roman{R}}$ 
to be a
representation? The first thing to notice is that, if 
$|i-j|\geq 2$, the
nontrivial parts of $\rho_{n,\roman{R}}(\sigma _i)$ and
$\rho_{n,\roman{R}}(\sigma _j)$  will not interfere with 
one another, so (4a)
is satisfied by construction, independently of the choice 
of R. As for (4b),
it is clear that we only need to look at the actions of 
$\roman{R} \otimes I_V$
and $I_V \otimes \roman{R}$ on $V^{\otimes 3}$, for if  $$ 
(\roman{R} \otimes
I_V)(I_V \otimes \roman{R})(\roman{R} \otimes I_V) =  (I_V 
\otimes
\roman{R})(\roman{R} \otimes I_V)(I_V \otimes \roman{R}), 
\leqno (7)$$ then
(4b) will be satisfied for {\it every \/} pair $\sigma_i, 
\sigma_j$ with
$|i-j| = 1$. Thus equation (7) is the clue to the 
construction. It is known as
the {\it quantum Yang-Baxter equation\/} (QYBE). Note that 
if R satisfies the
QYBE equation, then $\alpha \roman{R}$ will too, for every 
invertible element
$\alpha \in \Cal{E}$. By our earlier observations (see the 
discussion of
equations (5) and $(5^{\prime}))$, if the representation 
determined by R can be
used to define a link invariant, then the same will be 
true if R is replaced by
$\alpha$R. We will use this fact shortly, modifying our 
choice of R by
composing with an invertible scalar, whenever it 
simplifies the equations to do
so.

{\it Caution\/}. One must distinguish between (7), which 
we call the QYBE, and
a closely related equation $(7^*)$, which is also  
referred to as the QYBE. To
explain the latter and to show that both forms have the 
same geometric meaning,
let us return momentarily to braids. Our elementary braid 
generators could have
just as well been called $\sigma _{i,i+1}$ to stress the 
fact that the
nontrivial part of the braid involves the adjacent strands 
$i$ and $i+1$.
Relation (4b) in $B_3$ would then be  
$$\sigma _{12} \sigma _{23} \sigma _{12} = \sigma _{23} 
\sigma _{12}
\sigma _{12}.$$ 
An alternative way to number the braid generators would be 
to
color the strands with three colors, labeled $1,2,3$ and 
to let $\check{\sigma}_{i,j}$ denote a
positive crossing between the strands of color $i$ and 
color $j$, when those two colors are
adjacent in the braid projection. With that convention, we 
see from the bottom left
picture in Figure 14 on page 277 (which is the braid 
version 
of Reidemeister's move III) that (4b) would be
$$\check{\sigma}_{12} \check{\sigma}_{13} 
\check{\sigma}_{23} = \check{\sigma}_{23}
\check{\sigma}_{13} \check{\sigma}_{12}. \leqno (4.\roman 
b^*)$$  
Going back to operators acting on
$V^{\otimes{3}}$, let $\roman{R}_{12} = \roman{R} \otimes 
I_V$ and let $\roman{R}_{23} = I_V
\otimes \roman{R}.$  Then equation (7) may be rewritten  as 
$$\roman{R}_{12}\roman{R}_{23}\roman{R}_{12} = 
\roman{R}_{23}\roman{R}_{12}\roman{R}_{23}.$$
Now let $\check{\roman{R}}_{ij}$ be the image of 
$\roman{R}_{ij}$ under the automorphism of
$V^{\otimes 3}$ induced by the permutation of $\{ 1,2,3 
\}$ which maps $1$ to $i$ and $2$ to
$j$. The alternative form of the QYBE is 
$$\check{\roman{R}}_{12}\check{\roman{R}}_{13}\check{%
\roman{R}}_{23} =
\check{\roman{R}}_{23}\check{\roman{R}}_{13}\check{%
\roman{R}}_{12}.
\leqno (7^*)$$ In this
form, it occurs in nature in many ways, for example, in 
the theory of exactly solvable models in
statistical mechanics, where it appears as the {\it 
star-triangle 
relation \/} \cite {Bax}.

Not long after the discovery that the Jones polynomial 
generalized to the HOMFLY and Kauffman 
polynomials, workers in the area began to discover other, 
isolated cases of generalizations, 
all relating to existing known solutions to (7*).  Then a 
coherent theory emerged.
It turns out that solutions to $(7^*)$, and hence (7), can 
be constructed with the help of 
known solutions to yet a third equation---the {\it 
classical \/} Yang-Baxter equation
(CYBE). To explain that equation, let $\Cal{L}$ be a Lie 
algebra and let
$\roman{r} \in \Cal{L}^{\otimes 2}$. Let $\roman{r}_{12}$ 
be $\roman{r} \otimes 1 \in
\Cal{L}^{\otimes 3}$ and let $\roman{r}_{ij}$ be the image 
of $\roman{r}_{12}$ under the
automorphism of $\Cal{L}^{\otimes 3}$ induced by the 
permutation of $\{ 1,2,3 \}$ which sends
$1$ to $i$ and $2$ to $j$. Then r is said to be a solution 
to the CYBE if the
following holds: $$ \lbrack \roman{r}_{12}, \roman{r}_{23} 
\rbrack + \lbrack \roman{r}_{13},
\roman{r}_{23} \rbrack + \lbrack \roman{r}_{12}, 
\roman{r}_{13} \rbrack = 0. 
\leqno (7^{**})$$ 
The theory of the CYBE is well understood, its solutions 
having
been essentially classified. Two good references on the 
subject, both with 
extensive bibliographies, are 
\cite {Sem} and \cite {BD}. The possibility  of 
``quantization",
i.e., passage from solutions to the CYBE $(7^{**})$ to 
those for the QYBE $(7^*)$ was proved
by Drinfeld in a series 
of papers, starting with \cite {Dr1}. See \cite {Dr2} for 
a review
of the subject and another very extensive bibliography, 
including the definition of
a quantum group and a  development of the relevance of the 
theory of
quantum groups to this matter. That theory has much to do 
with the representations of 
simple Lie algebras. Explicit solutions to $(7^*)$ which 
are associated to the fundamental
representations of the nonexceptional classical Lie 
algebras may be found in the work of
Jimbo (\cite {Ji1} and especially \cite {Ji2}).  See \cite 
{WX} for a discussion of the
problem of finding ``classical'' solutions to the 
``quantum'' equation $(7^*)$.

Now, representations of $B_n$ do not always support a 
Markov trace. To obtain a
Markov trace from the representation  defined in (6), 
where R satisfies (7),
one needs additional data in the form of an {\it 
enhancement \/} of R  \cite
{Re1}. See Theorem 2.3.1 of [{Tu}]. The enhancement is a 
choice of invertible
elements $\mu _1,\dots,\mu _m \in \Cal{E}$ which determine 
a matrix  $\mu =
\roman{diag}(\mu _1,\dots,\mu _m)$ such that 
$$ \mu \otimes \mu \ \roman{commutes \ with}\ 
\roman{R} = \lbrack \roman{R}_{i_1i_2}^{j_1j_2} \rbrack, 
\tag 8a$$ 
$$ \sum_{j=1}^m
(\roman{R}^{\pm 1})_{i,j}^{k,j}\cdot\mu _j = \delta_{i,k} 
(\roman{Kroneker 
\ symbol).} \tag 8b$$
({\it A remark for the experts\/}. Theorem 2.3.1 of \cite 
{Tu} contains two factors,  $\alpha$ and
$\beta$. We have chosen them to be 1.
This possibility is discussed in Remark (i), \S3.3 of 
\cite {Tu}. We may do this because, as
noted earlier, if R is a solution to the QYBE, then 
$\alpha$R is too, for any invertible $\alpha
\in \Cal{E}$.)  

As it turns out, there are invertible elements 
$\mu_1,\dots,\mu_m \in \Cal{E}$
associated naturally to every solution to the QYBE which 
comes from quantum
groups and has the property that (8a) is true. This is 
proved in \cite {Ros,
Re 2, RT 1}. It is also known that if the $R$-matrix comes
from an irreducible representation of a quantum group, 
then (8b) is true with $
\delta_{i,k}$ replaced by $\alpha^{-1}\delta_{i,k}$ for 
some invertible element
$\alpha \in \scr E$.
Notice that (7) and (8a) still hold if R is
replaced by $\alpha$R. Thus, possibly replacing R by 
$\alpha$R, every solution
to the QYBE coming from quantum groups admits an 
enhancement.  

Finally, we describe how generalized Jones invariants are 
constructed out of the
data $(\roman{R}, n, \mu_1,\dots,\mu_m)$. Let $\rho = 
\rho_{n,\roman{R}}$. Let $\mu^{\otimes
n}$ denote the mapping $\mu \otimes \dots \otimes \mu : 
V^{\otimes n} \to V^{\otimes n}$. 
If $\beta \in
B_n$, define  $\TR(\rho(\beta))$ to be the ordinary matrix 
trace of $\rho(\beta) \cdot
\mu^{\otimes n}$. Notice that this implies that, if $1_n$ 
denotes the identity element in $B_n$,
then  $\TR(1_n)=(\mu_1+\cdots +\mu_m)^n$ 
is the matrix trace of $\mu^{\otimes n}$.  With this 
rescaling of the trace the
factor $z$ in equation (5) goes away (see the discussion 
following the statement of equation
(5)). The link invariant $F(\BK_{\beta})$ is given 
explicitly as $$ F(\BK_{\beta}) =
\TR(\rho(\beta) \cdot \mu^{\otimes n}) \.  \leqno (9)$$  
We leave it to the reader to verify that (8a) and (8b) 
ensure that
$F(\BK_{\beta})$ in (9) is invariant under Markov's two 
moves. 
(Alternatively, see Theorem 3.12 of \cite{Tu} for a 
proof.) The ring
$\Cal{E}$ (to the best of our knowledge) will always be a 
ring of Laurent
polynomials over the integers, the variable being some 
root of $q$. In the
special case of the HOMFLY polynomial (as defined in (3a) 
and (3b)), the
function $F(\BK_{\beta})$ of (9) turns out to depend  only 
on $q^{\pm 1}$ when
$\BK_{\beta}$ is a knot.  For links, square roots 
generally occur.  The 
invariant
$F(\bold K_\beta)$ in (9) is normalized so that $F(\bold 
O)=\mu_1+\cdots
+\mu_m$.

The HOMFLY and Kauffman polynomials, which are central to 
the subject, 
are related to the fundamental representations of the 
nonexceptional Lie
algebras of type $A^1_n$ and
$B^1_n, C^1_n, D^1_n, A^2_n$. Other knot and link 
polynomials, for example, 
the ones in \cite {WAD}, turn out to be the
HOMFLY and Kauffman polynomials of various cables on 
simple (noncabled) knots and links. Yet other 
completely new polynomials arise out of the 
representations of the 
remaining nonexceptional Lie
algebras and the exceptional Lie algebras. For example, 
see \cite {Ku} for the case of the exceptional
Lie algebra $G_2$ from the point of view 
of crossing-change formulas and \cite {Kal} for an 
interpretation of Kuperberg's polynomial via traces on 
representations of $B_n$. We understand that
the representations corresponding to the algebras $E_6, 
E_7, E_8$ have also been constructed, but we
do not have a reference. The important point is that it is 
possible to do so.

Summarizing our results in this section, we have described 
a very general construction
which leads to the construction of vast families of knot 
and link invariants. The proof, 
via Markov's theorem, that they do
not depend on the choice of the representative braid 
$\beta \in B_{\infty}$ 
is
straightforward and convincing. Unfortunately, however, 
neither the construction (which
depends upon the combinatorial properties needed to 
satisfy (7) or (7*) and 
(8)) nor the 
proofs we know of Markov's theorem (\cite {Be, Bi, Mo1}
and most recently  \cite {Mak})
give the
slightest hint as to the underlying topology; moreover, 
the constructions which
actually yield all possible solutions to (7) and (8) do 
not change that
situation at all. The crossing-change formulas in \S2 give 
even less insight.
Thus we can compute the simplest of the invariants by hand 
and quickly fill
pages with the coefficients and exponents of $J_q(\BK)$ 
for not-too-complicated
knots without having the slightest idea what they mean.

\heading 4. A space of all knots
\endheading

To proceed further, we need to change our point of view. 
The Alexander
polynomial, as we have seen, 
contains information about the topology of the
complement of a single knot. We now describe yet another 
method of constructing
knot invariants, discovered in 1989, which yields 
information about the
topology of an appropriately defined space of {\it all \/} 
knots. 

Following methods pioneered by Arnold (e.g., see the 
introduction to \cite
{Arn1}) and developed for the case of knots in 3-space by 
Vassiliev \cite {V},
we begin by introducing a change which is very natural in 
mathematics, shifting
our attention from the knot $K$, which is the image of 
$S^1$ under an embedding
$\phi:S^1 \rightarrow S^3$, to the embedding $\phi$ 
itself. A knot type $\BK$
thus becomes an equivalence class $\lbrace\phi\rbrace$ of 
embeddings of $S^1$ into
$S^3$. The space of all such equivalence classes of 
embeddings is disconnected,
with a component for each smooth knot type, and the next 
step is to enlarge it
to a connected space. With that is mind, we pass from 
embeddings to smooth
maps, thereby admitting maps  which have various types of 
singularities. Let
$\widetilde{\Cal{M}}$ be the space of all smooth maps from 
$S^1$ to $S^3$. This
space is connected and contains all knot types. Our space 
will remain
connected and will contain all knot types if we place two 
mild restrictions on
our maps. Let $\Cal{M}$ denote the collection of all $\phi 
 \in
\widetilde{\Cal{M}}$ such that $\phi(S^1)$ passes through 
a fixed point  $\ast$
and is tangent to a fixed direction at $\ast$.   The space 
$\Cal{M}$ has some
pleasant properties, the main one being that it can be 
approximated by certain
affine spaces, and these affine spaces contain 
representatives of all knot
types. The walls between distinct chambers in $\Cal{M}$ 
constitute the {\it
discriminant \/} $\Sigma$, i.e.,  
$$ \nomultlinegap\multline
\Sigma = \lbrace \phi \in \Cal{M}|\phi \text{
has a multiple point or a place where its derivative 
vanishes or}\\
\text{other singularities}\rbrace\.
\endmultline$$
The space $\Cal{M} - \Sigma$ is our
space of all knots.  ({\it Remark\/}. It is 
important to include in $\Sigma$ all maps
which have cusp singularities, for if not every knot could 
be changed to the
unknot by pulling it tight and, at the last moment, 
popping it!)

The group $\widetilde{H}^0(\Cal{M} - \Sigma ; \Bbb{Q})$ 
contains all
rational-valued knot-type invariants. (The ``tilde" over 
$H^0$ indicates that
we have normalized by requiring each invariant to take the 
value 0 on the
unknot.) Clearly, $\widetilde{H}^0(\Cal{M} -\Sigma; 
\Bbb{Q})$ contains enough
information to classify knots. For example, in principle 
it is possible to make
a list containing all knot types, ordered, say, by 
crossing number, and then to
delete repeats until knots which occur in distinct 
positions on the list have
distinct knot types. The position of a knot on the list 
would then define an
element in $\widetilde{H}^0(\Cal{M} - \Sigma; \Bbb{Q})$ 
which (by definition)
classifies knots.

Our space $\Cal{M}$ is much too big for direct analysis, 
but fortunately its very size makes
it extremely flexible and subject to approximation. Let 
$\Gamma^d \subset \Cal{M}$ be the
subspace of all maps from $\Bbb{R}^1$ to $\Bbb{R}^3$ given 
by
$$t \rightarrow (\phi _1(t), \phi_2(t),\phi_3(t)),$$ 
where each $\phi_i$ is a polynomial of the form 
$$t^{d+1} + a_{i1}t^d + \dots + a_{id}t$$
with $d$ even. As $t \rightarrow \pm \infty$,
the images $\phi(t)$ are seen to be asymptotic to
the direction $\pm(1,1,1)$, so the images are tangent to a 
fixed direction at $\ast =
\lbrace \infty \rbrace$. The space $\Gamma^d$ has three 
key properties:

(i) By the Weierstrass Approximation Theorem, every knot 
type occurs for some
$\phi \in \Gamma^d$ for  sufficiently large $d$. 

(ii) There are ways to embed $\Gamma^{d}$ in  
$\Gamma^{d'}, d <
d'$. For example \cite {V}, map $\Gamma^d \rightarrow 
\Gamma^{3d+2}$ by
reparametrizing $\phi(t)$ as $\phi(s^3 + s)$.   Thus, one 
may choose a sequence
of positive integers $d_1 < d_2 < d_3 < \cdots$ so that 
$$ \widetilde H^0(\Cal{M} - \Sigma; \Bbb{Q} ) =  
\lim_{\leftarrow} \widetilde
H^0 ( \Gamma^{d_i} - \Gamma^{d_i} \cap \Sigma; \Bbb{Q}). 
\leqno (10)$$

(iii) The fact that each $\phi \in \Gamma^d$ is determined 
by $3d$ real numbers
and that each
$3d$-tuple of real numbers determines a map in $\Gamma^d$ 
allows us to identify $\Gamma^d$ with
Euclidean space $\Bbb{R}^{3d}$. Alexander duality then 
applies, and 
$$ \widetilde H^0 ( \Gamma^{d_i} - \Gamma^{d_i} \cap 
\Sigma; \Bbb{Q}) \cong 
\bar H_{3d-1}(\Gamma^{d_i} \cap \Sigma; \Bbb{Q}).  \leqno 
(11) $$
where $\bar H_{3d-1}$ denotes reduced homology. In other 
words, the object of study is
no longer the topology of an individual component $\BK$ of 
$\Cal{M} - \Sigma$
 but instead
the topology of the walls $\Sigma$ which separate the 
components in our space 
of all knots. This change has an important consequence: it 
adds structure in
the form of the natural stratification of $\Sigma$ which 
occurs via the
hierarchy of its singularities. 

The discriminant $\Sigma$ is of course hideously 
complicated, since maps in
$\Sigma$ may have multiple points, tangencies, and all 
sorts of bad
singularities; but other studies in singularity theory 
have shown that certain
generic singularities predominate, and the next goal is to 
single these out.
Let $\Cal{M}_j \subset \Sigma$ be the subspace of maps 
$\phi \in \Sigma$ which
have $j$ transverse double points (and possibly other 
singularities too). Call
$\phi \in \Cal{M}_j$ a {\it j-embedding\/} if its only 
singularities are $j$
transverse double points. Let $\Sigma _j =  \lbrace \phi 
\in \Cal{M}_j/ \phi$
is not a $j$-embedding$\rbrace$. Then $\Sigma \supset 
\Sigma _1 \supset \Sigma
_2 \supset \cdots$.  Each map in $\Cal{M}_j - \Sigma_j$ 
has exactly $j$ generic
singularities, where a generic singularity is a transverse 
double point. 

The space $\Cal{M}_j - \Sigma_j$ is not connected. Its 
components are
singular knot types with exactly $j$ transverse double 
points. To illustrate how it might
enter our picture, let $\BK$ and $\BK^{\ast}$ be knot 
types which are neighbors in $\Cal{M} -
\Sigma$. Thus a single passage across $\Sigma$ suffices to 
change a representative $K$ of
$\BK$ to a representative $K^{\ast}$ of $\BK^{\ast}$. At 
the instant of passage across
$\Sigma$ (assuming that the passage is across a ``nice" 
part of $\Sigma$), one will obtain a
singular knot $K^1$ in $\Cal{M}_1$. This singular knot can 
be thought of as keeping
track of the crossing change. The singular knot $K^1$ will 
be in $\Cal{M}_1 - \Sigma _1$. 
The path from $\BK$ to $\BK^{\ast}$ will in general not be 
unique, so it may happen that two
distinct passages yield singular knots $K^1, N^1$ which 
are in distinct components of
$\Cal{M}_1 - \Sigma _1$. To pass from one to another, one 
must cross $\Sigma _1$, obtaining (at
the instant of passage) a singular knot $K^2 \in \Cal{M}_2 
- \Sigma_2$, and so forth.
Vassiliev studied this situation, applying the methods of 
spectral sequences to obtain
combinatorial conditions which determine certain families 
of invariants. These have
become known as {\it Vassiliev invariants \/}. They are 
stable values of elements in
the group  $\bar H_{3d-1}(\Gamma^d \cap \Sigma) \cong 
\widetilde H^0(\Gamma ^d - \Gamma ^d \cap
\Sigma)$, as $d \to \infty$, evaluated on a component 
$\BK$ of $\Cal{M} - \Sigma$. As such,
they are elements of $\widetilde H^0(\Cal{M} - \Sigma; 
\Bbb{Q})$. 

A Vassiliev invariant $v_i(\BK)$ of $order\ \ i$ \ takes 
into account
information about the singular knot types near $\BK$ in 
$\Cal{M} - \Sigma, \
\Cal{M}_1 - \Sigma _1,\ \Cal{M}_2 - \Sigma 
_2,\dots,\Cal{M}_i - \Sigma _i$,
where $i$ \ is a positive integer. The sum of two $i$\<th 
order invariants is
another invariant of order $i$. In fact, these invariants 
lie in a
finite-dimensional vector space $V_i$ of dimension $d_i$. 
Clearly, these
invariants have much to do with crossing changes.

The manuscript \cite {V} ends with a combinatorial recipe 
for the calculation
of $v_i(\BK)$.  
   
\heading 5. Axioms and initial data for Vassiliev invariants
\endheading

In 1990, when preprints of \cite {V} were first 
circulated, topologists were in
possession of an overflowing cornucopia of knot and link 
invariants. In
addition to the Jones polynomial and its generalizations, 
we mention the knot
group invariants of \cite {Wa}, the energy invariants of 
\cite {FH}, and the
algebraic geometry invariants of \cite {CCG}, all of which 
seemed to come
from mutually unrelated directions! In addition, there 
were generalizations of
the Jones invariants to knotted graphs \cite {Y1} and, 
finally, the numerical
knot invariants of \cite {V}.  

The first hint that the Vassiliev and Jones invariants 
might be related was
that both extended to invariants of singular knots. A 
second hint was a tiny
bit of evidence: the first nontrivial Vassiliev invariant 
is $v_2$, which was
identified in  \cite {V} as coinciding with the second 
coefficient in the
Conway-Alexander polynomial. In fact, that coefficient has 
another
interpretation, as the second derivative of the Jones 
polynomial evaluated at
1.  Studying \cite {V} and having in mind a possible 
relationship with Jones
invariants, the author and Lin were led to reformulate the 
results in
\cite {V} in terms of a set of ``axioms and initial data". 
We
describe them in this section.

Let $v: \Cal{M} - \Sigma \to \Bbb{Q}$ denote a function 
from the space $\Cal{M}
- \Sigma$ of all knots to the rational numbers. This 
function determines a
Vassiliev invariant if it satisfies the following 
conditions: 

\subheading{Axioms}
 To state the first axiom, we note that the singular knots 
in \cite {V}
 are subject to an equivalence relation which is known as 
{\it rigid vertex
isotopy \/}---a neighborhood of each double point of the 
singular knot spans an open disc
$\Bbb{R}^2 \subset \Bbb{R}^3$---and this disc must be 
preserved under admissible isotopies. For
example, if one begins with the standard diagram of the 
trefoil and of its mirror image (both shown
in Figure 4) and then replaces any one crossing in each 
diagram by a double point, one will obtain
two singular knots which are equivalent under isotopy but 
not under rigid vertex
isotopy. 
This restriction is quite natural if one thinks of 
singular knots as objects which
keep track of crossing changes in a passage from one knot 
to a neighbor in $\Cal{M} - \Sigma$.
Knowing that there is such a disc, or plane, it then makes 
sense to speak of the two {\it
resolutions \/} (see Figure 5)  at a singular point $p$ of 
the singular knot $K^j$: We denote
these two resolutions $K_{p_+}^{j-1}$ and $K_{p_-}^{j-1}$ 
respectively.  
The first axiom is
a type of crossing-change formula:  $$ v(K^j_p) =
v(K^{j-1}_{p_+}) - v(K^{j-1}_{p_-}). \tag 10a$$   
This axiom in itself places no
constraints on $v$ as a knot invariant. Rather, it serves 
to define $v$ on singular knots,
assuming that it is known on knots.

The second axiom is what makes $v$ computable: 
$$ \exists i
\in Z^+\ \text{such\ that\ }v(\bold K^j) = 0 \text{ 
if}\ j> i.    \tag 10b$$
The smallest such $i$ is the {\it orde \/r} of $v$. To 
stress it, we call our invariant
$v_i$ from now on.      

\subheading{Initial data}
In addition to (10a) and (10b) we need initial data. The 
first
piece of initial data relates to the normalization 
mentioned earlier:
$$v_i(\bold{O}) = 0   \.\tag  10c$$
To describe the second piece of initial data, we need a 
definition. A singular point on a
knot diagram will be called {\it nugatory \/} if its 
positive and negative resolutions define
the same knot type, in the obvious manner indicated in 
Figure 10. It is clear that if we are to
obtain a true knot invariant, its value on $v_i(\BK 
^{j-1}_{p_+})$ and $v_i(\BK ^{j-1}_{p_-})$
must agree when $p$ is a nugatory crossing, so by (10a) 
the initial 
data must satisfy 
$$v_i(\BK ^j_p) = 0 \text{ if } p \text{ is a nugatory 
crossing}. \tag  10d$$

\topspace{16pc}
\caption{{\smc Figure} 9. {\rm Actuality table for $i=2$.}}

\np The final piece of initial data is in the form of a 
table, but before we
can describe it we need to discuss a point that we avoided 
earlier: The space
$\Cal{M}_j - \Sigma_j$ has a natural decomposition into 
components, such that
two singular knots cannot define the same singular knot 
type if they belong to
distinct components. To make this assertion precise, let 
$\BK ^j$ be a singular
knot of order $j$, i.e., the image of $S^1$ under a 
$j$-embedding $\phi \in
\Cal{M}_j - \Sigma_j$. Then $\phi^{-1}(\BK ^j)$ is a 
circle with $2j$
distinguished points, arranged in pairs, where two 
distinguished points are
paired if they are mapped to the same double point on $\BK 
^j$. The $[j]$-{\it
configuration \/} which $\BK ^j$ {\it respects \/} is the 
cyclically ordered
collection of point pairs. We will use a picture to define 
it, i.e., a circle
with arcs joining the paired points, as in the top row of 
Figure 9, where we
show the two possible [2]-configurations, together with a 
choice of a singular
knot which respect each. The initial data must take 
account of the following
(see \cite {ST 1} for a proof):

\proclaim{Lemma 1}  Two singular knots $K^j_1, K^j_2$ 
become equivalent after an appropriate series of
crossing changes if and only if they respect the same 
$[j]$-configuration.
\endproclaim  

\fighere{5pc}
\caption{{\smc Figure 10}. {\rm Nugatory crossing.}} 

The table which we now construct to complete the initial 
data is called an {\it
actuality table \/}.  Figure 9 is an example, when $i=2$. 
It gives the values
of $v_i(\BK ^j)$ for a representative collection of 
singular knots of order $j
\leq i$. The table contains a choice of a singular knot 
$K^j$ which respects
each $[j]$-configuration, for $j=1,2,\dots,i$.  The choice 
is arbitrary (however,
the work in completing the table will be increased if poor 
choices are made).
Next to each $K^j$ in the table is the configuration it 
respects, and below it
is the value of $v_i(\BK ^j)$. These values are, of 
course, far from arbitrary,
and the heart of the work in \cite {V} is the discovery of 
a finite set of
rules which suffice to determine them. The rules turn out 
to be in the form of
a system of linear equations. The unknowns are the values 
of the functionals on
the finite set of singular knots in the actuality table. 
The linear equations
which hold between these unknowns are consequences of the 
local equations
(which may be thought of as crossing-change formulas) 
illustrated in Figure 11.
These equations are not difficult to understand: use (10a) 
to resolve each
double point into a sum of two crossing points. Then each 
local picture in
Figure 11 will be replaced by a linear combination of four 
pictures. The
equations in Figure 11 will then be seen to reduce to a 
sequence of
applications of Reidemeister's third move. See \S3 of 
\cite {BL} for a
description of the method that allows one to write down 
the full set of
equations. See \S 2.4 of \cite {BN} for a proof that 
solutions to the
equations, in the special case when $j=i$, may be 
constructed out of 
information about the irreducible representations of 
simple Lie algebras.  It
is not known whether the methods of \cite {BN} yield {\it 
all\/} solutions in the
case $j=i$.  The extension of the solutions for the case 
$j=i$ to the cases
$2\le j\le i-1$ must be handled by the less routine 
methods described in \cite
{BL}, at this time. 

An example should suffice to illustrate that (10a)--(10d) 
and the actuality table allow one to
compute $v_i(\BK)$ on all knots. For our example we 
compute $v_2(\BK)$ when $\BK$ is the
trefoil knot. The\vadjust{\fighere{14.5pc}
\caption{{\smc Figure 11}. {\rm Crossing change formula 
for Vassiliev invariants of
related singular knots.}}}
first
picture in Figure 12 shows our representative of the 
trefoil, with a
crossing which is marked. Changing it, we will obtain the 
unknot $\bold{O}$.

\topspace{19.5pc}
\caption{{\smc Figure 12}. {\rm Computing $v_2(\BK)$ for 
the trefoil knot.}}

The crossing we selected is positive  and so (10a) yields 
$$v_2(K) =v_2(O) + v_2(N^1)$$ 
where $N^1$ is the indicated singular knot. It does
not have the same singular knot type as the singular knot 
$K^1$ in the table, so (using the
lemma) we introduce another crossing
 change to modify it to the singular knot in the
table which respects the unique [1]-configuration. In so 
doing, we obtain a singular knot $N^2$
with two singular points. It does not have the same 
singular knot type as the representative in
the table which respects the same [2]-configuration, but 
by (10b) that does not matter. Thus,
using (10b), the calculation comes to an end in finitely 
many steps. 

\heading 6. Singular braids
\endheading

In \S 3 we showed that any generalized Jones invariant may 
be obtained from a Markov
trace on an appropriate family of finite-dimensional 
matrix representations of the braid
groups. Up to now, braids have not entered the picture as 
regards Vassiliev invariants,
but that is easy to rectify. To do so, we need to extend 
the usual notions of braids and
closed braids to singular braids and closed singular 
braids. 

A representative $K ^j$ of a singular knot or link $\BK 
^j$ will be said to be
a {\it closed singular braid \/} if there is an axis $A$ 
in $\Bbb{R}^3$ (think
of it as  the $z$-axis) such that if $K ^j$ is 
parametrized by cylindrical
coordinates $(z, \theta)$ relative to $A$, the polar angle 
function restricted
to $K ^j$ is monotonic increasing. This implies that $K^j$ 
meets each
half-plane $\theta = \theta{_0}$ in exactly $n$ points, 
for some $n$. In
\cite {Al1} Alexander proved the well-known\ fact that 
every knot or link $\BK$ may be
so represented, and we begin our work by extending his 
theorem to singular
knots and links. 

\proclaim{Lemma 2} Let $K^j$ be an arbitrary 
representative of a singular knot or link
$\BK ^j$. Choose any copy $A$ of $\Bbb{R}^1$ in $\Bbb{R}^3 
- K^j$. Then $K^j$ may be deformed to a
closed singular n-braid, for some n, with $A$ \ as 
axis.\endproclaim

\demo{Proof} Regard $A$ as the $z$-axis in $\Bbb{R}^3$.  
After an isotopy of
$K^j$ in \ $\Bbb{R}^3 - A$ \ we may assume that $K^j$ is 
defined by a diagram
in the $(r, \theta)$-plane.  By a further isotopy we may 
also arrange that each
singular point $p_k$ has a neighborhood $N(p_k) \in K^j$ 
such that the polar
angle function restricted to $\bigcup_{k=1}^{j} N(p_k)$ is 
monotonic
increasing. The proof then proceeds exactly as in  \cite 
{Al1}, vis: Modify \
$K^j - \bigcup_{k=1}^{j} N(p_k)$ \ to a piecewise linear 
family of arcs
$\Cal{A}$, subdividing the collection if necessary so that 
each $\alpha \in
\Cal{A}$ contains at most one undercrossing or 
overcrossing of the knot
diagram. After a small isotopy we may assume that the 
polar angle function is
nonconstant on each $\alpha \in \Cal{A}$. Call an arc 
$\alpha \in \Cal{A}$ {\it
bad \/} or {\it good \/} accordingly as the polar angle 
function is increasing
or decreasing on $\alpha$. If there  are no bad arcs, we 
will have a closed
braid, so we may assume there is at least one, say 
$\beta$. Modify $K^j$ by
replacing $\beta$ by two good edges $\beta _1 \cup \beta 
_2$ as in Figure 13.
The only possible obstruction is if the interior of the 
triangle which is
bounded by $\beta \cup \beta{_1} \cup \beta{_2}$  is 
pierced by the rest of
$K^j$, but that may always be avoided by choosing the new 
vertex $\beta{_1}
\cap \beta{_2}$ so that it lies very far above (resp.\  
below) the rest of
$K^j$ if the arc $\beta$ contains an under- (resp.\ over-) 
crossing of the
diagram. Induction on the number of bad edges completes 
the proof.  
\qed\enddemo 

In view of Lemma 2, every singular knot in the actuality 
table of \S5 may be
chosen to be a closed singular braid. To continue, we 
split these closed
singular braids open along a plane $\theta = \theta{_0}$ 
to ``open" singular
braids, which we now define. In \S3 we described a 
geometric braid as a pattern
of $n$ interwoven strings in $\Bbb{R}^2 \times I \subset 
\Bbb{R}^3$ which join
$n$ points, labeled $1,2,\dots,n$ in  $\Bbb{R}^2 \times 
\{0\}$ to the
corresponding points in $\Bbb{R}^2 \times \{1\}$, 
intersecting each 
intermediate plane $\Bbb{R}^2 \times \{t\}$ in exactly $n$ 
 points. To extend
to singular braids, it is only necessary to weaken the 
last condition to allow
finitely many values of $t$ at which the braid meets 
$\Bbb{R}^2 \times \{t\}$
in $n-1$ points instead of $n$ points. Two singular braids 
are equivalent if
they are isotopic through a sequence of singular  braids, 
the isotopy fixing
the initial and end points of each singular braid strand. 
Singular braids
compose like ordinary braids: concatenate two patterns,  
erase the middle
plane, and rescale. 

\midspace{7pc}
\caption{{\smc Figure 13}. {\rm Replacing a bad arc by two 
good arcs.}}

Choose any representative of an element of $SB_n$. After 
an isotopy we may
assume that distinct double points occur at distinct 
$t$-levels. From this it
follows that $SB_n$ is generated by the elementary braids
$\sigma_1,\dots,\sigma_{n-1}$ and the elementary singular 
braids $\tau_1,\dots,
\tau_{n-1}$ of Figure 8. We distinguish between the 
$\sigma_i$\<'s and the
$\tau_i$\<'s by calling them {\it crossing  points \/} and 
{\it double points
\/} respectively in the singular braid diagram. Both 
determine double points in
the projection. 

The manuscript \cite {Bae} lists defining relations in 
$SB_n$ as:
$$[\sigma_i, \sigma_j] = [\sigma_i, \tau_j] = [\tau_i, 
\tau_j] = 0\quad
\text{if }  |i-j| \geq 2,
\tag 11a$$
$$[\sigma_i, \tau_i] = 0,  \tag  11b$$
$$ \sigma_i\sigma_j\sigma_i = \sigma_j\sigma_i\sigma_j 
\quad \text{if }  |i-j| =
1,  \tag 11c$$
$$ \sigma_i\sigma_j\tau_i = \tau_j\sigma_i\sigma_j \quad 
\text{if } |i-j| = 1, 
\tag 11d$$
where in all cases $1 \leq i, j \leq n-1.$ The same set of 
relations also
occurs 
in \cite {KV} as generalized Reidemeister moves.The 
validity of these relations is easily
established via pictures; for example,
see Figure 14 for special cases of (11a)--(11d).
To the best of our knowledge, 
however, there is not even a sketch of a proof that they 
suffice
in the
literature, so we sketch one now, as it will be important 
for us that no additional relations are
needed. 

\proclaim{Lemma 3} The monoid $SB_n$ is generated by
$\{\sigma_i,\tau_i ; \ 1 \leq i \leq n-1\}$. Defining 
relations are
{\rm (11a)--(11d)}.
\endproclaim

\demo{Proof} We have already indicated a proof that the 
$\sigma_i$\<'s  and the
$\tau_i$\<'s generate $SB_n$, so the only question is 
whether every relation is
a consequence of (11a)--(11d). 

We regard our braids as being defined by diagrams. Let 
$\overline{z},
\overline{z}'$ be singular braids which represent the same 
element of $SB_n$,
and let $\{ \overline{z}_s; \ s \in I \}$ denote the 
family of singular braids
which join them. The fact that the intermediate braid 
diagrams $\overline{z}_s$
have no\vadjust{\fighere{14.5pc}
\caption{{\smc Figure 14}. {\rm Relations in $SB_n$.}}}
triple points implies that
there is a well-defined order of the
singularities along each braid strand which is preserved 
during the isotopy,
and this allows us to set up a 1-1 correspondence between 
double points in
$\overline{z}$ and $\overline{z}'$. We now examine the 
other changes which
occur during the isotopy. Divide the $s$-interval [0,1] 
into small
subintervals, during which exactly one of the following 
changes occurs in the
sequence of braid diagrams:

(i) Two double points in the braid projection interchange 
their $t$-levels. See
relations (11a) and (11b) and Figure 14.
  
(ii) A triple point in the projection is created 
momentarily as a ``free"
strand crosses a double point or a crossing point in the 
projection. See (11d)
and Figure 14.

(iii) New crossing points in the knot diagram are created 
or destroyed. See
Reidemeister's second move in Figure 2. 

\np All possible cases of (i) are described by relations 
(11a) and (11b).
Noting that the $\sigma_i$\<'s are invertible and that the 
mirror image of
$\sigma_i$ (resp.\ $\sigma_i^{-1}, \tau_i)$ is 
$\sigma_i^{-1}$ (resp.\
$\sigma_i, \tau_i)$, it is a simple exercise to see that 
consequences of (11c)
and (11d) cover all possible cases of (ii). As for (iii), 
if we restrict to
small $s$-intervals about the instant of creation or 
destruction, these will
occur in pairs and be described by the trivial relation 
$\sigma_i \sigma_i^{-1}
= \sigma_i^{-1} \sigma_i = 1$. Outside of these special  
$s$-intervals the
singular braid diagram will be modified by isotopy, which 
contributes no new
relations. Thus relations (11a)--(11d) are defining 
relations for $SB_n$.
\qed\enddemo

Now something really interesting happens. Let 
$\widetilde{\sigma}_i$ denote the
image of the elementary braid $\sigma_i$ under the natural 
map from the braid
group $B_n$ to its group algebra $\Bbb{C}B_n$. 

\proclaim{Theorem 2 {\rm (cf.\ \cite {Bae, Lemma 1})}\rm}
 The map $\eta : SB_n \to
\Bbb{C}B_n$ which is defined by $\eta (\sigma_i) = 
\widetilde{\sigma}_i, \eta 
(\tau_i) = \widetilde{\sigma}_i
- \widetilde{\sigma}_i^{-1}$ is a monoid homomorphism.
\endproclaim 

\demo{Proof} Check to see that relations (11a)--(11d) are 
consequences of the
braid relations in $\Bbb{C}B_n$.   
\qed\enddemo

\proclaim{Corollary 1} Every finite-dimensional matrix 
representation $\rho _n:
B_n \to \GL _n(\Cal{E})$ extends to a representation 
$\tilde{\rho}_n$ of
$SB_n$, defined by setting $\tilde{\rho}_n(\tau_i) = 
\rho_n(\sigma_i) -
\rho_n(\sigma_i^{-1}).$ 
\endproclaim

\demo{Proof} 
Clear.
\qed\enddemo

\np Thus, in particular, all R-matrix representations of 
$B_n$ extend to
representations of the singular braid monoid $SB_n$.

\rem{Remark \RM1} Recall that in Figure 11 we gave picture 
examples of some of
the relations which need to be satisfied by a functional 
on the knot space in
order for the indices in an actuality table to determine a 
knot type invariant.
Unlike the relations which we depicted in Figure 14, the 
ones in Figure 11 are,
initially, somewhat mysterious. However, if one passes via 
 Theorem 2 to the
group algebra $\Bbb{C}B_n$ of the braid group, replacing 
each $\tau _i$ by
$\sigma _i - \sigma _i^{-1}$, it will be seen that these 
relations actually
hold in the algebra, not simply in the space of 
VBL-functionals. This fact is
additional evidence of the naturality of the map $\eta$ of 
Theorem 2 and indeed
of the Vassiliev construction. We conjecture that the 
kernel of $\eta$ is
trivial, i.e., that a nontrivial singular braid in the 
monoid $SB_n$ never maps
to zero in the group algebra $\Bbb{C}B_n$.\endrem

\rem{Remark \RM2} Various investigators, for example,  
Kauffman in \cite
{Kau1}, have considered a somewhat different monoid which 
we shall call the
{\it tangle  monoid \/}, obtained by adding the 
``elementary tangle"
$\varepsilon _i$ of Figure 8 to the braid group. The 
tangle monoid, however,
does not map homomorphically to $\Bbb{C}B_n$, and  one 
must pass to quotients
of $\Bbb{C}B_n$, for example,  the so-called Birman-Wenzl 
algebra \cite {BW,
We1}, to give it algebraic meaning. In that sense the 
tangle monoid appears to
be less fundamental than the singular braid monoid.\endrem

\heading 7. The proof of Theorem 1
\endheading

We are now ready to prove Theorem 1. We refer the reader 
to the introduction
for its statement.

Theorem 1 was first proved by the author and Lin in \cite 
{BL} for the special
cases of the HOMFLY and Kauffman polynomials and then in 
full generality for
all quantum group invariants in \cite {L}. The proof we 
give here is new. It is
modeled on the proof for the special cases in \cite {BL}. 
We like it because it
is simple and because it gives us an opportunity to show 
that the braid groups,
which are central to the study of the Jones invariants, 
are equally useful in
Vassiliev's setting. The tools in our proof are the 
R-matrix representations of
\S 3, the axioms and initial data of \S 5, and the 
singular braid monoid of \S
6. Theorem 1 also is implicit in \cite {Bae}, which was 
written simultaneously
with and independently of this manuscript. The techniques 
used there are very
similar to ours, but the goal is different. 

\demo{Proof of Theorem \RM1} By hypothesis, we are given a 
quantum group
invariant $\Cal{J}_q: \Cal{M}-\Sigma \to \Cal{E}$, where 
$\Cal{E}$ denotes a
ring of Laurent polynomials over the integers in powers of 
$q$ (or in certain
cases powers of roots of $q$). Using Lemma 2, we find a  
closed braid
representative $K_{\beta}$ of $\BK$, $\beta \in B_n$. We 
then pass to the
R-matrix representation $\rho _{n,R}: B_n \to \GL 
_{m_n}(\Cal{E})$ associated
to $\Cal{J}_q$. By Corollary 1 that representation extends 
to a representation
$\tilde{\rho}_{n,R}$ of $SB_n$.  By equation (9) the 
Laurent polynomial
$\Cal{J}_q(\BK)$ is the trace of $\rho_{n,R}(\beta) 
\cdot\mu$, where $\mu$ is
the enhancement of $\rho_{n,R}$.\enddemo

As was discussed in \S3, our representation $\rho_{n,R}$ 
is determined by the
choice of a matrix $R$ which acts on the  vector space 
$V^{\otimes 2}$. Lin
notes in  Lemma 1.3 of \cite {L} that on setting $q=1$ the 
matrix $R =
(R^{j_1j_2}_{i_1i_2})$  goes over to $\ID_V \otimes 
\ID_V$.  From this it
follows that $\rho_{n,R}(\sigma_i)$ (which acts on 
$V^{\otimes n}$) has order 2
at $q=1$. Hence, we conclude that $\rho_{n,R}$ goes over 
to a representation of
the symmetric group if we set $q=1$, with $\sigma_i$ going 
to the transposition
$(i,i+1)$. In particular, this means that at $q=1$ the 
images under
$\rho_{n,R}$ of $\sigma_i$ and $\sigma_i^{-1}$ will be 
identical, which, in
turn, means that the images of $\sigma_i$ and 
$\sigma_i^{-1}$ under
$\rho_{n,R}$ coincide at $q=1$.

Armed with this knowledge, we change variables, as in the 
statement of Theorem 1,
replacing $q$ by $e^x$. Expanding the powers of $e^x$ in 
its Taylor series, the image of an
arbitrary element $\beta$ under $\rho_{n,R}$ will be a 
matrix power series
$$\rho_{n,R}(\beta) = M_0(\beta)
+ M_1(\beta) + M_2(\beta) + \cdots    \leqno (12)$$
where each $M_i(\beta) \in \GL _{m_n}(Q)$. 

\proclaim{Lemma 5 {\rm (cf.\ \cite {Bae, Corollary 1})}\rm}
In the extended representation,
$M_0(\tau_i) = 0.$ \endproclaim

\demo{Proof} Since $M_0(\sigma_i) = M_0(\sigma_i^{-1})$, 
the assertion follows.
\qed\enddemo

Now let us  turn our attention to the power series 
expansion of $\Cal{J}_q$,
i.e., to $\Cal{J}_x(\BK) = \sum_{i=0}^{\infty} u_i(\BK) 
x^i.$  The coefficients
$u_i(\BK)$ in this series, as $\BK$ ranges over all knot 
types, determine a
functional $u_i:\Cal{M} - \Sigma \to Q$. We wish to prove 
that $u_i$ is a
Vassiliev invariant of order $i$. 

By \S5 it suffices to prove that, if we use (10a) to 
extend the definition of
$u_i$ to singular knots, then (10b)--(10d) will be 
satisfied and a consistent
actuality table exists.  The first thing to notice is 
that, since we began with
a knot-type invariant $\Cal{J}_q$, the functional $u_i$ is 
also a knot-type
invariant. From this it follows that its extension to 
singular knots is also
well defined, so using our knowledge of $u_i$ on knots we 
can fill in the
actuality table. 

The second observation is that (10c) is  satisfied, 
because every generalized
Jones invariant satisfies $\Cal{J}_q(\bold{O}) \equiv 1$, 
and from this it
follows that $\Cal{P}_x(\bold{O}) \equiv 1$. As for (10d), 
it is also
satisfied, for if not $u_i$ could not be a knot-type 
invariant. Thus the only
problem which remains is to prove that $u_i$ satisfies 
axiom (10b). However,
notice that by Lemma 5 we have
$$\tilde{\rho}_{n,R}(\tau_i) = 
M_1(\tau_i)x + M_2(\tau_i)x^2 + M_3(\tau_i)x^3 +\cdots,$$
and from this it follows that if $\BK^j$ is a singular 
knot which has $j$
singular points, then any singular closed braid 
$K_{\gamma}^j, \ \gamma \in
SB_r$, which represents $\BK^j$ will also have $j$ 
singularities.  The singular
braid word $\gamma$ will then contain $j$ elementary 
singular braids. From this
it follows that
$$\tilde{\rho}_{n,R}(\gamma) =
M_{j}(\gamma )x^j + M_{j+1}(\gamma )x^{j+1} + \cdots.$$ 
The coefficient of $x^i$ in this power series is $u_i(\bold
 K^j)$. But then, $u_i(\bold K^j) = 0$ 
if $i<j$, i.e., axiom (10b) is also satisfied and the 
proof of Theorem 1 is complete.  
\qed 

\heading 8. Open problems
\endheading

Our story is almost ended, so it is appropriate to 
recapitulate and ask what we
have gained by our interpretation of the quantum group 
invariants as Vassiliev
invariants? The goal has been to bring topology into the 
picture, and indeed we
have done so because we have shown that, when the Jones 
invariant of a knot
$\BK$ is expanded in a power series in $x$, the 
coefficient of $x^i$ gives
information about certain stable homology groups of the 
discriminant $\Sigma$,
in a neighborhood of $\BK$. The information concerns its 
structure at ``depth"
$i$; that is, of course, only a beginning. The 
discriminant $\Sigma$ is a
complicated subset of an infinite-dimensional space,  
$\Cal{M}$, and there
seems to be  no way to begin to visualize $\Sigma$. 
Indeed, the work has just
begun.

The study of Vassiliev invariants is fairly new. We now 
review some of the
things that have been learned about them during the past 
few years. A natural
place to start is with the special cases when $i$ is 
small, and this was
already done in \cite {V}. The invariant $v_1(\BK)$ is 
zero for all knots
$\BK$. The first nontrivial Vassiliev invariant has order 
2, and there is a
one-dimensional vector space of such invariants. However, 
$v_2(\BK)$ was
well known to topologists in other guises before either 
Jones or Vassiliev came
on the scene, vis:

(i) It is the second coefficient in Conway's version of 
the Alexander
polynomial \cite {C}.

(ii) Its mod 2 reduction is the Arf invariant of a knot, 
which has to do with
cobordism.

\np It now has other interpretations too:

(iii) It is the ``total twist" of a knot, as defined in 
\cite {LM}. It may be
computed from  the following recursive formula: 
$$ v_2(K_{p_+}) - v_2(K_{p_-}) = \Cal{L}k(K_{p_0}),$$ 
where $\Cal{L}k$ denotes the linking number of the
two-component link $K_{p_0}$. 

(iv) It is the second derivative of the Jones polynomial 
$J_q(\BK)$, evaluated
at $q=1$. 

Unfortunately, all of these facts do not allow us to 
understand the
significance of $v_2$ for the topology of $\Sigma$. There 
is much work to be
done. As for $v_3(\BK)$, it also belongs to a 
one-dimensional vector space of
invariants, but to the best of our knowledge no one has 
yet identified it as a
classical invariant---indeed, nothing much seems to be 
known about it. 

With regard to the higher-order invariants, recall that 
Vassiliev invariants of
order $i$ belong to a linear vector space $V_i$. This 
space is the space of
$i$\<th order invariants modulo those of order $i-1$. The 
first characteristic
of this vector space to question
 is its dimension $m_i$, for $i > 2$. This seems to be a
deep and difficult combinatorial problem for arbitrary i, 
and at this writing
the best we can do is to compute, the issue being the 
construction of the
actuality table for an invariant of order $i$. This 
problem divides naturally
into two parts:  the first is to determine the Vassiliev 
invariants of the
singular knots $K^i$ in the top row. By axiom (10b) they 
only depend upon 
the $[i]$-configuration which they respect, so their 
determination is easier than
the corresponding problem for the remaining rows. The 
former problem is solved
by setting up and solving the system of linear equations 
(3.11) of \cite {BL}.
The dimension of the space of solutions is the sought-for 
integer $m_i$. One
year ago Bar Natan wrote a computer program which listed 
the distinct
$[i]$-configurations and calculated the dimensions (top 
row only) for $i \leq 7$.
Very recently Stanford developed a different computer 
program which checked
that all of Bar Natan's ``top row" solutions, for $i \leq 
7$, extend to the
remaining rows of the actuality table, i.e., to solutions 
to the somewhat more
complicated set of equations (3.17) of \cite {BL}, so Bar 
Natan's numbers
are actually the dimensions $m_i$ that we seek. The 
results of the two 
calculations are:

\midinsert
\noindent$$
\table\nobox
i & 1 & 2 & 3 & 4 & 5 & 6 & 7\\
m_i & 0 & 1 & 1 & 3 & 4 & 9 & 14\endtable$$
\endinsert

The data we just gave leads us to an important question. 
We proved, in Theorem
1, that  
$$\{\text{Quantum\ group\ invariants}\} \subseteq    \{
\text{Vassiliev\ invariants}\}.$$ 
The question is: Is the inclusion proper?  
The data is relevant because it was thought that the 
question might be answered
by showing that for some fixed $i$ the quantum group 
invariants spanned a
vector space of dimension $d_i < m_i$; however, a  
dimension count to $i=6$
shows that there are enough linearly independent 
invariants coming from quantum
groups to span the vector spaces $V_i$. Bar Natan's 
calculations showed that
$d_7$ is at least 12, whereas $m_7 = 14$; however, the 
data on the quantum group
invariants is imprecise because the invariants which come 
from the
nonexceptional Lie algebras begin to make their presence 
felt as $i$
increases. The only one of those which has been 
investigated, to date,
is $G_2$.  

For $i=8$ the computations themselves create difficulties. 
Bar Natan's
calculation is close to the edge of what one can do, 
because to determine $m_8$
he found he had to solve a linear system consisting of 
334,908 equations in
41,874 unknowns (the number of distinct [8]-configurations 
with no separating
arcs). His approximate calculation shows that the  
solution space is of
dimension 27. However, even if his answer is correct, we 
would still need to
compute the rest of the actuality table before we could be 
sure that $m_8 = 27$
rather than $m_8 \leq 27$. 

Vassiliev invariants actually form an algebra, not simply 
a sequence of  vector
spaces, because  the product of a Vassiliev invariant of 
order $p$ and one of
$q$ is a Vassiliev invariant  of order $p+q$. This was 
proved by Lin
(unpublished), using straightforward methods, and more 
indirectly by Bar Natan
\cite {BN}. Thus the dimension $\hat{m}_i$ of new 
invariants is in general
smaller than $m_i$, because the data in our table includes 
 invariants which
are products of ones of smaller order. It is a simple 
matter to correct the
given data to find $\hat {m}_i$. For example, an invariant 
of order $4$ could
be the product of two invariants of order $2$; so when we 
correct for the fact
that there is a one-dimensional space of invariants of 
order $2$, we see that
$\hat{m}_4 = 3 - 1 = 2$.  For $i=1,2,3,4,5,6,7$ we find 
that $\hat{m}_i =
0,1,1,2,3,5,8$, i.e., the beginning of the Fibbonaci 
sequence! This caused some
excitement until Bar Natan's computation of $m_8$ showed 
that $\hat{m}_8$ was
at most 12, not 13. The asymptotic behavior of $m_i$ as $i 
\to \infty$ is a
very  interesting problem indeed.

We can approach the question of whether the inclusion is 
proper from a
different point of view. One of the earliest problems in 
knot theory concerned
the fundamental symmetries which are always present in the 
definition of knot
type. We defined knot type to be the topological 
equivalence class of the pair
$(S^3, K)$ under homeomorphisms which preserve the given 
orientation on both
$S^3$ and $K$. A knot type is called {\it amphicheiral \/} 
if it is equivalent
to the knot type obtained by reversing the orientation of 
$S^3$ (but not $K$)
and  {\it invertible \/}  if it is equivalent to the knot 
type obtained by
reversing the orientation on $K$ (but not $S^3$). As noted 
earlier, Max Dehn
proved in 1913 that nonamphicheiral knots exist \cite 
{De}, but remarkably, it
took over forty years before it became known that 
noninvertible knots exist \cite
{Tr}. The relevance of this matter to our question is: 
While the quantum group
invariants detect nonamphicheirality of knots, they do not 
detect
noninvertibility of knots. So, if we could prove that 
Vassiliev invariants
distinguished a single noninvertible knot from its 
inverse, the answer to our
question, Is the inclusion proper? would be  yes. We know 
we cannot 
answer the question this way for $i \leq 7$, and as noted 
above $i=8$ presents
serious computational difficulties. On the other hand, the 
theoretical problem
seems to be unexpectedly subtle. Thus, at this moment, the 
matter of whether
Vassiliev invariants ever detect noninvertibility remains 
open.

Setting aside empirical evidence and unsolved problems, we 
can ask some easy
questions which will allow us to sharpen the question of 
whether Jones or
Vassiliev invariants determine knot type. As was noted in 
\S2, there are three
very intuitive invariants which have, to date, proved to 
be elusive: the
crossing number $c(\BK)$, the unknotting number $u(\BK)$, 
and the
braid index
$s(\BK)$. Clearly these determine functionals on the space 
$\Cal{M}- \Sigma$,
and so they determine elements in the group $\widetilde 
{H}^0(\Cal{M}-
\Sigma)$. Vassiliev invariants lie in a sequence of 
approximations to
$\widetilde {H}^0(\Cal{M}- \Sigma)$. So a reasonable 
question to ask is, Are
$c(\BK), \ u(\BK)$, and $s(\BK)$ Vassiliev invariants? 
Theorem 5.1 of \cite
{BL} shows that $u(\BK)$ is not, and the proof given there 
is easily modified
to show that $c(\BK)$ and $s(\BK)$ are not either. So, at 
the very least, we
have learned that there are integer-valued functionals on 
Vassiliev's space of
all knots which are not Vassiliev invariants. This leaves 
open the question of
whether there are sequences of Vassiliev invariants which 
converge to these
invariants.

Another question which has been asked is, How powerful are 
the Vassiliev
invariants, if we restrict our attention to invariants of 
bounded order? The
answer to that question is, not very good, based on 
examples which were
discovered, simultaneously and independently, by Lin \cite 
{Li2} and by 
Stanford \cite {St2}. We now describe Stanford's 
construction, which is
particularly interesting from our point of view because it 
is based on  the
closed-braid approach to knots and links. See Remark (ii) 
for a description of
Lin's construction.

To state his theorem, we return to braids. Let $P_n$ be 
the pure braid group,
i.e., the kernel of the natural homomorphism from $B_n$ to 
the symmetric group
$S_n$. \ The groups of the lower central series $\{ 
P_n^k;\  k=1,2,\dots\}$ of
$P_n$ are defined inductively by $P_n^1 = P_n,\  P^k_n = 
\lbrack P_n, P_n^{k-1}
\rbrack$. Notice that if $\beta \in B_n$, with the closed 
braid
$\widehat{\beta}$ a knot, then  $\widehat{\alpha  \beta}$ 
\ will also be a knot
for every $\alpha \in P_n$.  

\proclaim{Theorem 3 \cite {St2}\rm}  Let $\BK$ be any knot 
type, and let
$K_{\beta}$, where $\beta \in B_n$, be any closed braid 
representative of
$\BK$.  Choose any $\alpha \in P_n^k$. Then the Vassiliev 
invariants of order
$\leq k$ of the knots $K_{\beta}$ and $K_{\alpha \beta}$ 
coincide. 
\endproclaim

\rem{Remarks}  (i) Using the results in \cite {BM1}, 
Stanford has constructed
sequences $\alpha_1, \alpha_2$, $\dots$ of 3-braids such that 
the knot types
obtained via Theorem 3 are all distinct and prime. 
Intuition suggests that
distinct $\alpha _j$\<'s will always give distinct knot 
types, but at this
writing that has not been proved. 

(ii) One may choose $K_{\beta}$ in Theorem 3 to represent 
the unknot and
obtain infinitely many distinct knots all of whose 
Vassiliev invariants of
order $\leq k$, for any $k$, are zero. Lin's construction 
gives other such
examples. In particular, he proves that, if $K$ is any 
knot and if $K(m)$ is
its $m$\<th iterated (untwisted) Whitehead double, then 
all Vassiliev
invariants of order $\leq (m+1)$ of $K(m)$ are zero. This 
allows him to
construct infinitely many composite knots which have the 
same properties as
Stanford's knots. It is not clear whether his construction 
is a special case of
Stanford's construction. 

(iii) The special case of the unknot and the one-variable 
Jones polynomial is
interesting. Theorem 3 says that if we expand the Jones 
polynomials of the
knots in the collection in power series and truncate those 
series by cutting
off all terms of $x$-degree bigger than $k$, the result 
will be zero. This is
far from saying that the polynomials themselves are 
trivial. Indeed, the
question of whether there is a nontrivial knot which has 
the same one-variable
Jones polynomial as the unknot remains an open problem and 
an intriguing one.
Note that our lack of knowledge about this problem is in 
striking contrast to
the control mathematicians now have over the Alexander 
polynomial:
understanding its topological meaning, we also know 
precisely how to construct
knots with Alexander polynomial  1.

(iv) For a brief period (before a mistake was discovered 
in the proof) it
seemed as if there might be an affirmative answer to the 
following question:
Given any two distinct knots $\BK$ and $\BK^*$, does there 
exist a sequence of
Vassiliev invariants $\{ w_i\ |\ i = 1,2,3,\dots\}$ and an 
integer $N$ such
that $w_i(\BK) \not= w_i(\BK^*) \ \forall i \geq N$?  
Notice that even if the
answer was yes, it would not solve the knot problem, 
because if we had
explicit examples $\BK$ and $\BK^*$ and if we knew the 
sequence, we still would
not know how large $N$ had to be. So after letting the 
computer run all weekend
without a definitive answer, we would not know whether the 
knots were really the
same or whether we simply had given up too soon. However, 
this is probably the
very best that one could hope for from the algebraic 
invariants.

(v) By contrast to all of this, the work of Haken \cite 
{Ha} and the work of
Hemion \cite {He} show that there is an algorithm which 
distinguishes knots.
In recent years efforts have been directed at making that 
algorithm workable
(for example,  see \cite {JT} for a discussion of recent 
results), 
but much work 
remains to be done before it could be considered 
``practical", even for the
simplest examples.

(vi) The work in \cite {BM2} and related papers referenced 
therein is aimed at
a different algorithm which is based upon the theory of 
braids. At present it
has resulted in a very rapid hand calculation for 
definitively distinguishing
knots and links which are defined by closed braids with at 
most 3 strands, with
partial results for the general case. 

(vii) Among the special cases for which an effective 
classification scheme
exists we mention (in addition to the work just cited on 
links defined by
closed 3-braids) the  cases of links with 2-bridges and of 
algebraic (in the
sense of algebraic geometry) links. For more information 
on these and other
classical topics, see \cite {Ro1} and \cite {BZ}.

(viii) An extensive guide to the ``pre-Jones" literature 
may be found in \cite
{BZ}.\endrem
   
\heading Acknowledgments \endheading 
Many people helped us to understand this subject. We 
single out, especially: V.
I. Arnold (who first stimulated our interest in Vassiliev 
invariants),
Xiao-Song Lin, and Dror Bar Natan. Discussions with 
Vaughan Jones, Vladimir
Turaev, Michael Polyak, Alex Sossinsky, and Ted Stanford 
have also been very
helpful to us. We also thank Oleg Viro for telling us 
about the parallel theory
of knot invariants which was developed by M. Gusarov \cite 
{Gu}.



\Refs
\widestnumber\key{Koh 2}

\ref
\key     Al1          
\by	    J. W. Alexander            
\paper  A lemma on systems of knotted curves         
\jour  Proc. Nat. Acad. Sci. U.S.A.          
\vol    9          
\yr   1923                                 
\pages  93--95                       
\endref



\ref		
\key   Al2          
\by J. W. Alexander          
\paper  Topological invariants of knots and links         
\jour   Trans. Amer. Math. Soc.          
\vol     30          
\yr     1928                                 
\pages  275--306                       
\endref


\ref
\key      Arn1          
\by       V. I. Arnold (ed.)
\book   Theory of singularities and its applications
\publ  Amer. Math. Soc.
\publaddr Providence, RI
\bookinfo Adv. in Soviet Math.,
 vol. 1
\yr    1990                                   
\endref

\ref		
\key   Arn2          
\bysame 
\paper  The cohomology ring of the colored braid group     
\jour   Notes Acad. Sci. USSR          
\vol     5          
\yr     1969                                 
\pages  227--232                       
\endref


\ref
\key   Art          
\by   E. Artin            
\paper  Theorie der Z\"opfe
\jour   Hamburg Abh.          
\vol    4          
\yr      1925                                 
\pages      47--72                       
\endref


\ref
\key   Bae          
\by    J. Baez            
\paper   Link invariants of finite type and perturbation 
theory         
\jour   Lett. Math. Phys.          
\toappear                      
 \endref

\ref
\key   Bax          
\by   R. J. Baxter            
\book  Exactly solvable models in statistical mechanics
\publ Academic Press
\publaddr London
\yr 1982         
\endref


\ref
\key  BD            
\by  A. A. Belavin and V. G. Drinfeld               
\paper  On classical Yang-Baxter equation for simple Lie 
algebras         
\jour  Funct. Anal. Appl.
\vol 16          
\yr 1982                                
\pages 1--10                       
\endref

\ref
\key    Be         
\by    D. Bennequin            
\paper  Entrlacements et \'equations de Pfaffe         
\jour  Asterisque          
\vol  107--108          
\yr   1983                                 
\pages 87--161                       
\endref



\ref
\key          Bi 
\by           J. S. Birman
\paper Braids, links and mapping class groups
\publ         Princeton Univ. Press
\jour         Ann. of Math Stud., vol. 84
\publaddr     Princeton, NJ, 1974
\endref

\ref
\key      BL         
\by    J. S. Birman and X. S. Lin            
\paper  Knot polynomials and Vassiliev invariants          
\jour    Invent. Math
\toappear                       
\endref

\ref
\key BM1              
\by  J. S. Birman and W. W. Menasco               
\paper Studying links via closed braids {\rm III:} 
Classifying links 
which are closed $3$-braids          
\jour  Pacific J. Math.              
\toappear                           
\endref

\ref
\key BM2              
\bysame               
\paper Studying links via closed braids {\rm VI:} A 
non-finiteness theorem          
\jour  Pacific J. Math.,
 vol. 156
\yr 1992
\pages 265--285
\endref


\ref
\key BN           
\by   D. Bar Natan             
\paper  On the Vassiliev knot invariants         
\paperinfo   Harvard Univ., preprint 1992
\endref



\ref
\key    BW          
\by    J. S. Birman and H. Wenzl            
\paper  Braids, link polynomials and a new algebra         
\jour   Trans. Amer. Math. Soc.          
\vol    313          
\yr    1989                                 
\pages    249--273                       
\endref

\ref
\key    BZ          
\by    G. Burde and H. Zieschang            
\book  Knots         
\publ   de Gruyter          
\publaddr  Berlin and New York         
\yr    1985                       
\endref




\ref
\key     C          
\by  J. H. Conway            
\paper  An enumeration of knots and links and some of 
their algebraic properties         
\inbook  Computational Problems in Abstract Algebra        
\ed  J. Leech
\publ  Pergamon Press
\publaddr New York 
\yr  1970                                 
\pages  329--358                       
\endref

\ref
\key   CCG          
\by   D. Cooper, M. Culler, H. Gillet, D. D. Long, and P. 
B. Shalen            
\paper  Plane curves associated to character varieties of 
$3$-manifolds         
\paperinfo Univ. of California, Santa Barbara, preprint 1991
\endref


\ref
\key   CGL
\by   M. Culler, C. McA. Gordon, J. Leucke and P. B. 
Shalen            
\paper  Dehn surgery on knots         
\jour  Ann. of Math. (2)          
\vol  125          
\yr  1987                                
\pages  237--300                       
\endref

\ref
\key De             
\by  M. Dehn                
\paper Die beiden Kleeblattschlingen         
\jour Math. Ann. 
\vol  75          
\yr 1914                                
\pages 402--413                       
\endref


\ref
\key  Dr1          
\by  V. G. Drinfeld            
\paper  On a constant quasi-classical solution to the 
quantum 
Yang-Baxter equation         
\jour  Doklady Acad. Nauk SSSR          
\vol  273          
\yr  1983                                 
\pages  531--550                       
\endref


\ref
\key  Dr2          
\bysame            
\paper  Quantum groups         
\inbook  Proc. ICM Berkeley
\publ Amer. Math. Soc.
\publaddr Providence, RI          
\vol  1          
\yr 1987                                 
\pages  798--820                       
\endref

\ref
\key FH             
\by  M. Freedman and Z. H. He                
\paper On the energy of knots and unknots         
\paperinfo  Univ. of California, San Diego, preprint 1992  
\endref


\ref
\key  FHL
\by  P. Freyd, J. Hoste, W. B. R. Lickorish, K. 
Millett, A. Ocneanu, and D. Yetter            
\paper  A new polynomial invariant of knots and links      
\jour  Bull. Amer. Math. Soc. 
\vol  12          
\yr  1985                                 
\pages  239--246                       
\endref



\ref
\key  FW          
\by  J. Franks and R. F. Williams           
\paper  Braids and the Jones polynomial         
\jour  Trans. Amer. Math. Soc. 
\vol  303          
\yr  1987                                 
\pages  97--108                       
\endref

\ref
\key  GL          
\by  C. McA. Gordon and J. Leucke          
\paper Knots are determined by their complements        
\jour  Bull. Amer. Math. Soc. 
\vol  20 
\yr   1989                                
\pages  83--88                       
\endref


  
\ref
\key Go          
\by  L. Goeritz            
\paper Bermerkungen zur knotentheorie         
\jour  Hamburg Abh.          
\vol  10          
\yr  1934                                 
\pages  201--210                       
\endref

\ref
\key Gu
\by M. Gusarov
\paper A new form of the Conway-Jones 
polynomial for knots via von Neumann algebras
\jour  Zap. Nauchn. Sem. Leningrad Otdel. Mat. Inst. 
Steklov. (LOMI)
\vol  193
\yr 1991
\pages 4--9
\endref


\ref
\key  Ha              
\by  W. Haken               
\paper \"Uber das Hom\"omorphieproblem der 
$3$-Mannigfaltigkeiten          
\jour Math. Z. 
\vol 80
\yr  1962
\pages  89--120                           
\endref

\ref
\key  He              
\by  G. Hemion               
\paper On the classification of homeomorphisms of 
$2$-manifolds and the classification of $3$-manifolds      
\jour Acta Math.              
\vol  142
\yr 1979
\pages  123--155                           
\endref


\ref
\key   Ho         
\by  J. Hoste          
\paper A polynomial invariant of knots and links         
\jour Pacific J. Math.          
\vol 124 
\yr 1986                                 
\pages  295--320                       
\endref

\ref
\key   Je         
\by  L. Jeffrey          
\paper Chern-Simons-Witten invariants of lens spaces and 
torus bundles,
and the semiclassical approximation         
\paperinfo Oxford Univ., preprint 1992
\endref


\ref
\key Ji1          
\by  M. Jimbo            
\paper  A $q$-difference analogue of $U(g)$ and the 
Yang-Baxter equation         
\jour Lett. Math. Phys. 
\vol  10          
\yr  1985                                 
\pages  63--69                       
\endref


\ref
\key   Ji2         
\bysame            
\paper Quantum $R$-matrix for the generalized Toda system  
\jour Comm. Math. Physics          
\vol  102          
\yr  1986                                 
\pages  537--547                       
\endref


\ref
\key  Jo1         
\by  V. F. R. Jones            
\paper  Index for subfactors         
\jour  Invent. Math.          
\vol  72          
\yr  1983                                 
\pages  1--25                       
\endref


\ref
\key  Jo2          
\bysame            
\paper  Braid groups, Hecke algebras and type 
$\roman{II}_1$ factors         
\inbook  Geometric Methods in Operator Theory          
\bookinfo Pitman Res. Notes Math. Ser. 
\publ Longman Sci. Tech., Harlow 
\vol 123          
\yr 1986                                 
\pages  242--273                       
\endref



\ref
\key   Jo3         
\bysame            
\paper A polynomial invariant for knots via von Neumann 
algebras         
\jour Bull. Amer. Math. Soc. 
\vol  12          
\yr  1985                                 
\pages  103--111                         
\endref


\ref
\key  Jo4          
\bysame            
\paper Hecke algebra representations of braid groups and 
link polynomials         
\jour  Ann. of Math. (2)
\vol  126          
\yr  1987                                 
\pages  335--388                       
\endref

\ref
\key Jo5              
\bysame               
\paper On knot invariants related to some statistical 
mechanics models          
\jour Pacific J. Math.              
\vol 137 
\yr 1989
\pages  311--334                           
\endref


\ref
\key JR
\by V. Jones and M. Rosso
\paper Invariants of torus knots derived from quantum groups
\jour   Abstracts Amer. Math. Soc., Abstract no. 874-16-124
\yr  1992
\endref

\ref
\key    JT          
\by    W. Jaco and J. L. Tollefson            
\paper  Algorithms for the complete decomposition of a 
closed $3$-manifold         
\paperinfo    Univ. of Conn. at Storres, preprint 1992
\endref


\ref
\key  Kal          
\by E. Kalfagianni            
\paper The G$_2$ link invariant         
\paperinfo  Columbia Univ., preprint 1992
\endref


\ref
\key  Kan          
\by  T. Kanenobu            
\paper  Infinitely many knots with the same polynomial     
\jour  Proc. Amer. Math. Soc.
\vol  97 
\yr  1986                                 
\pages  158--162                       
\endref


\ref
\key  Kau 1          
\by L. H. Kauffman            
\paper  An invariant of regular isotopy         
\jour Trans. Amer. Math. Soc. 
\vol  318 
\yr   1990                                
\pages  417--471                       
\endref


\ref
\key Kau 2          
\bysame            
\paper States models and the Jones polynomial         
\jour   Topology          
\vol  26 
\yr  1987                                 
\pages 395--407                       
\endref

\ref
\key KM          
\by R. Kirby and P. Melvin            
\paper The $3$-manifold invariants of Witten and 
Reshetikhin-Turaev
for $SL(2,\Bbb{C})$         
\jour Invent. Math.          
\vol  105          
\yr  1991                                 
\pages  473--545                       
\endref


\ref
\key  Koh1         
\by  T. Kohno          
\paper Monodromy representations of braid groups and 
Yang-Baxter equations 
\jour Ann. Inst. Fourier (Grenoble)
\vol 37 
\yr 1987
\pages 139--160                      
\endref

\ref
\key Koh 2         
\bysame          
\paper Linear representations of braid groups and classical 
Yang-Baxter equations        
\inbook Contemp. Math., vol. 78
\publ Amer. Math. Soc.
\publaddr Providence, RI
\yr 1988
\pages 339--363                      
\endref

\ref
\key Kon         
\by M. Kontsevich
\paper Vassiliev's knot invariants          
\paperinfo Adv. in Soviet Math.
\toappear
\endref



\ref
\key  Ku          
\by  G. Kuperberg            
\paper  The $1,0,1,1,4,10$ Ansatz         
\paperinfo preprint
\endref

\ref
\key  KV          
\by  L. Kauffman and P. Vogel            
\paper  Link polynomials and a graphical calculus         
\jour  J. Knot Theory Ramifications          
\vol 1  
\yr  1992                                
\pages  59--104                       
\endref


\ref
\key La         
\by  R. Lawrence            
\paper  Homological representations of the Hecke algebras  
\jour  Comm. Math. Phys. 
\vol  135          
\yr  1990                                 
\pages  141--191                       
\endref


\ref
\key  Li1          
\by   X. S. Lin           
\paper  Vertex models, quantum groups and Vassiliev's knot 
invariants         
\paperinfo Colum\-bia Univ.,  preprint 1991 
\endref

\ref
\key Li2              
\bysame               
\paper Finite type link invariants of $3$-manifolds        
\paperinfo Columbia Univ.,  preprint 1992
\endref


\ref
\key   LM         
\by  W. B. R. Lickorish and K. Millett            
\paper  A polynomial invariant for knots and links         
\jour Topology          
\vol  26          
\yr  1987                                
\pages  107--141                       
\endref


\ref
\key Mak         
\by  G. S. Makanin           
\paper  On an analogue of the Alexander-Markov theorem     
\jour  Math. USSR Izv. 
\vol  34  
\yr  1990                                
\pages  201--211                      
\endref


\ref
\key  Mar         
\by  A. A. Markov           
\paper \"Uber die freie \"Aquivalenz geschlossener
Z\"opfe         
\jour  Recueil Math. Moscou          
\vol  1          
\yr  1935                                 
\pages  73--78                       
\endref

\ref
\key Me              
\by W. W. Menasco               
\paper A proof of the Bennequin-Milnor unknotting 
conjecture          
\paperinfo SUNY at Buffalo, preprint 1992
\endref


\ref
\key  Mo1         
\by   H. Morton           
\paper  Threading knot diagrams         
\jour  Math. Proc. Cambridge Philos. Soc.          
\vol 99          
\yr 1986                                
\pages  247--260                       
\endref


\ref
\key Mo2          
\bysame            
\paper Seifert circles and knot polynomials         
\jour Math. Proc. Cambridge Philos. Soc.           
\vol  99          
\yr  1986                                 
\pages  107--109                       
\endref






\ref
\key  Mu          
\by  K. Murasugi            
\paper  The Jones polynomial and classical conjectures in 
knot theory         
\jour  Topology          
\vol  26          
\yr  1987                                 
\pages  187--194                       
\endref

\ref
\key  Pi          
\by  S. Piunikhin            
\paper  Weights of Feynman diagrams and Vassiliev knot 
invariants          
\paperinfo 
Moscow State Univ., preprint 1992
\endref

\ref
\key Pr              
\by    C. Procesi               
\paper  The invariant theory of n by n matrices          
\jour  Adv. in Math.              
\vol  19
\yr  1976
\pages  306--381                           
\endref


\ref
\key  PT          
\by   J. Przytycki and P. Traczyk         
\paper  Conway algebras and skein equivalence of links     
\jour   Proc. Amer. Math. Soc. 
\vol  100          
\yr  1987                                 
\pages  744--748                       
\endref

\ref
\key  Re1          
\by   N. Yu. Reshetikhin         
\paper Quasitriangular Hopf algebras and invariants of 
tangles          
\jour  Leningrad Math. J.         
\vol  1 
\yr  1990                                 
\pages  491--513                       
\endref

\ref
\key  Re2          
\bysame        
\paper Quantized universal enveloping algebras, the 
Yang-Baxter 
equation and invariants of links          
\jour  LOMI preprints         
\vol  E-4-87, E-17-87
\yr 1988                       
\endref



\ref
\key Rol
\by D. Rolfsen
\book Knots and links         
\publ Publish or Perish        
\publaddr Berkeley, CA 
\yr 1976          
\endref

\ref
\key Ros              
\by M. Rosso               
\paper Groupes quantique et mod\`eles \`a  vertex 
de V. Jones en th\'eorie des noeuds          
\jour C. R. Acad. Sci. Paris               S\'er. I Math. 
\vol 307 
\yr  1988
\pages 207--210                           
\endref

\ref
\key RT          
\by N. Reshetikhin and V. G. Turaev            
\paper  Invariants of $3$-manifolds via link polynomials 
and quantum groups         
\jour  Invent. Math.          
\vol  103          
\yr  1991                                 
\pages  547--597                       
\endref

\ref
\key Sei             
\by H. Seifert                
\paper Verschlingungsinvarianten         
\jour Sber. Preuss. Akad. Wiss. 
\vol 26          
\yr 1933                                
\pages 811--828                       
\endref

\ref
\key Sem             
\by M. Semenov Tian-Shansky                
\paper What is a classical R-matrix\,{\rm ?}
\jour  Funct. Anal.  Appl. 
\vol 17          
\yr 1983                                
\pages  259--270                       
\endref

\ref
\key           Si
\by            J. Simon
\paper         How many knots have the same group\/{\rm ?}
\jour          Proc. Amer. Math. Soc. 
\vol           80 
\yr            1980                     
\pages         162--166               
\endref

\ref
\key  St1          
\by  T. Stanford            
\paper  Finite type invariants of knots, links and graphs  
\paperinfo  Columbia Univ., preprint 1992 
\endref

\ref
\key St2          
\bysame            
\paper  Braid commutators and Vassiliev invariants         
\paperinfo   Columbia Univ., preprint 1992                 
\endref

\ref
\key  Ta          
\by  P. G. Tait            
\paper On knots {\rm I, II} and {\rm III }
\inbook Scientific Papers of P. G. Tait          
\vol 1
\publ Cambridge Univ. Press          
\publaddr Cambridge and New York 
\yr 1988                                  
\pages  273--347                       
\endref

\ref
\key Tr             
\by H. Trotter                
\paper Non invertible knots exist         
\jour Topology 
\vol  2          
\yr  1964                                
\pages  275--280                       
\endref

\ref
\key  Tu          
\by   V. Turaev           
\paper  The Yang-Baxter equation and invariants of links   
\jour  Invent. Math.          
\vol  92          
\yr  1988                                 
\pages   527--553                       
\endref

\ref
\key         V1
\by           V. A. Vassiliev 
\paper        Cohomology of knot spaces
\inbook       Theory of Singularities and its Applications 
(V. I. Arnold, ed.)
\pages        23--69
\publ         Amer. Math. Soc. 
\publaddr     Providence, RI
\yr           1990
\endref

\ref
\key  V2        
\bysame            
\paper Topology of complements to discriminants and loop 
spaces         
\inbook Theory of Singularities and its Applications
\ed V. I. Arnold
\publ Amer. Math. Soc.          
\publaddr  Providence, RI          
\yr  1990                                 
\pages  9--21                       
\endref

\ref
\key Wa          
\by  M. Wada            
\paper  Group invariants of links         
\jour Topology          
\vol 31 
\yr 1992
\pages 399--406                       
\endref

\ref
\key WAD          
\by  M. Wadati, Y. Akutsu, and T. Deguchi            
\paper Link polynomials and exactly solvable models         
\inbook Non-Linear Physics           (Gu, Li, and Tu, eds.)
\publ Springer-Verlag
\publaddr Berlin          
\yr 1990                       
\endref

\ref
\key We1          
\by  H. Wenzl            
\paper  Quantum groups and subfactors of type B, C and D   
\jour Comm. Math. Phys. 
\vol  133   
\yr  1990                                 
\pages  383--432                       
\endref

\ref
\key  We2          
\bysame            
\paper Representations of braid groups and the quantum 
Yang-Baxter equation         
\jour  Pacific J. Math.          
\vol  145 
\yr  1990                                 
\pages  153--180                       
\endref

\ref
\key Wi          
\by E. Witten            
\paper  Quantum field theory and the Jones polynomial      
\jour  Comm. Math. Phys.           
\vol  121          
\yr  1989                                 
\pages  351--399                       
\endref

\ref
\key Wh          
\by W. Whitten            
\paper  Knot complements and groups         
\jour  Topology          
\vol  26 
\yr  1987                                 
\pages  41--44                       
\endref

\ref
\key WX          
\by A Weinstein and P. Xu            
\paper  Classical solutions of the quantum Yang-Baxter 
equation         
\jour  Comm. Math. Phys.          
\vol  148          
\yr  1992                                
\pages  309--343                       
\endref

\ref
\key Y1          
\by S. Yamada            
\paper  An invariant of spacial graphs         
\jour  J. Graph Theory          
\vol  13          
\yr  1989                                
\pages  537--551                       
\endref

\ref
\key Y2          
\bysame            
\paper  The minimum number of Seifert circles equals the 
braid index         
\jour  Invent. Math.           
\vol  89          
\yr  1987                                 
\pages  346--356                       
\endref
\endRefs
\enddocument